\documentclass[11pt, a4paper]{article}
\usepackage{amsmath}
\usepackage{amsfonts}
\usepackage{amssymb}
\usepackage{ulem}

\usepackage{amsthm}
\usepackage{url}

\usepackage{cite}
\usepackage{fancyhdr}

\usepackage{graphicx, color}
\usepackage{hyperref}

\begin{document}

\title{Convex Pentagon Tilings and Heptiamonds, I}
\author{ Teruhisa SUGIMOTO$^{ 1), 2)}$ and Yoshiaki ARAKI$^{2)}$ }
\date{}
\maketitle

{\footnotesize

\begin{center}
$^{1)}$ The Interdisciplinary Institute of Science, Technology and Art

$^{2)}$ Japan Tessellation Design Association

Sugimoto E-mail: ismsugi@gmail.com, Araki E-mail: 
yoshiaki.araki@tessellation.jp
\end{center}

}

\medskip

{\small
\begin{abstract}
\noindent
In 1995, Marjorie Rice discovered an interesting tiling 
by using convex pentagons. The authors discovered novel properties of the 
tilings of convex pentagon which Rice used in the discovery. As a result, 
many new convex pentagon tilings (tessellations) were found. The convex 
pentagon tilings are related to tilings by heptiamonds.
\end{abstract}


}

\section{Introduction}
\label{section1}

\begin{figure}[htbp]
 \centering\includegraphics[width=14.5cm,clip]{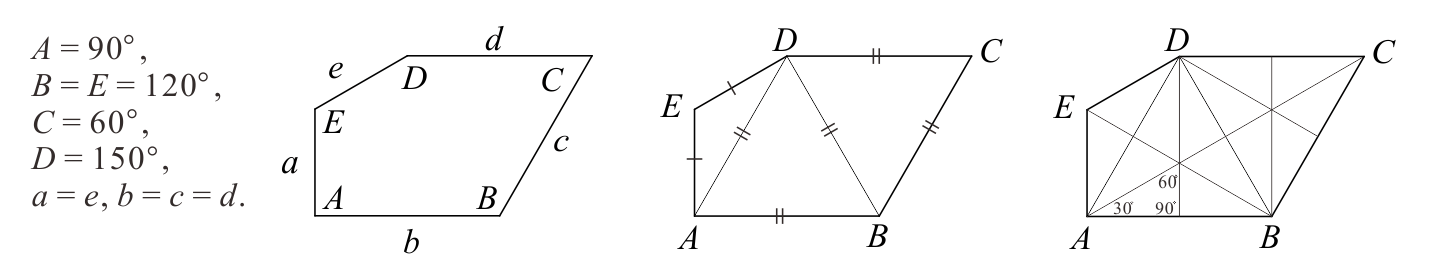} 
  \caption{{\small 
TH-pentagon that belongs to both Type 1 and Type 5.} 
\label{fig1}
}
\end{figure}

\begin{figure}[htbp]
 \centering\includegraphics[width=15cm,clip]{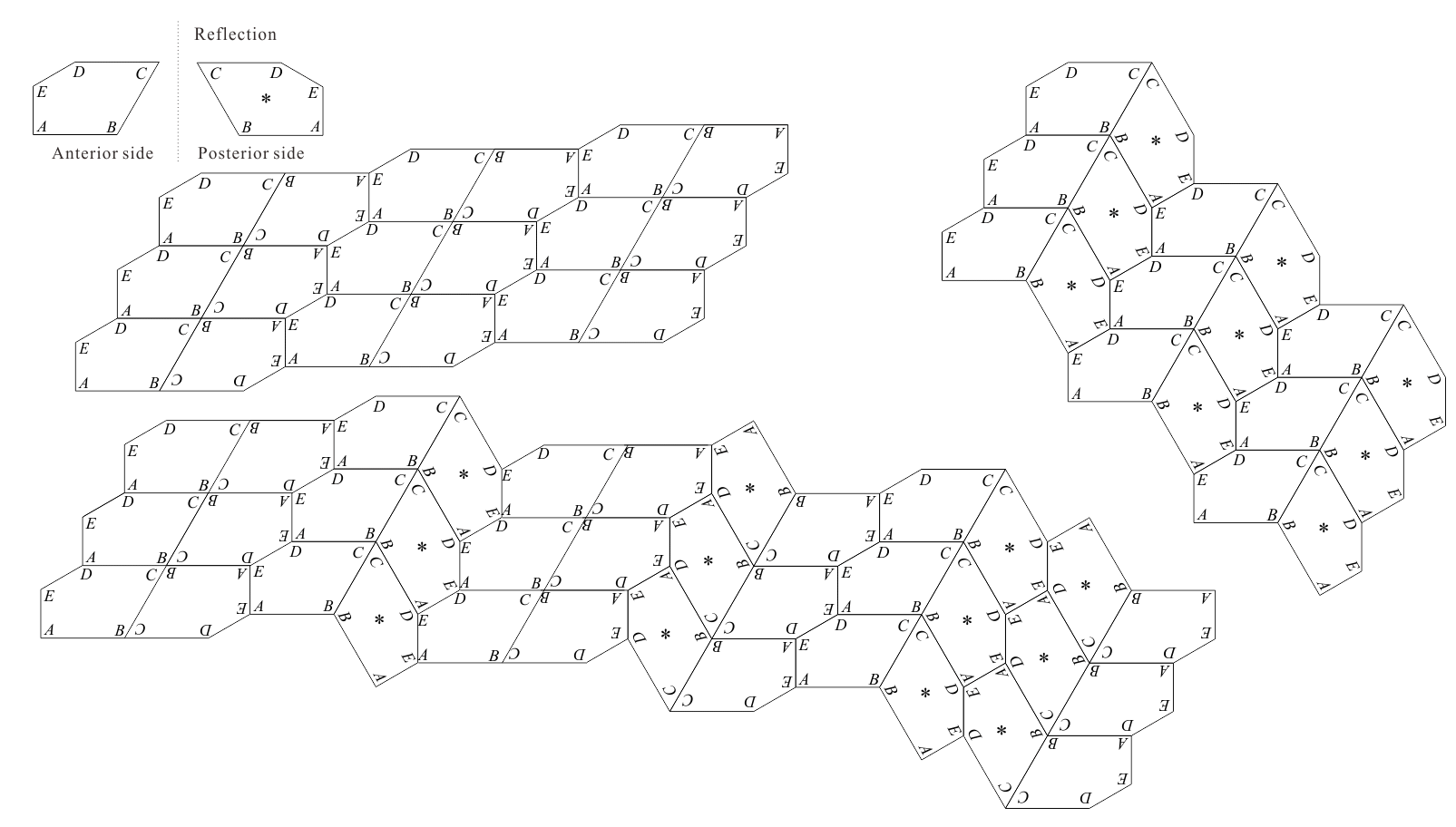} 
  \caption{{\small 
Representative tiling of Type 1 and variations of Type 1 tilings.} 
\label{fig2}
}
\end{figure}

\begin{figure}[htbp]
 \centering\includegraphics[width=15cm,clip]{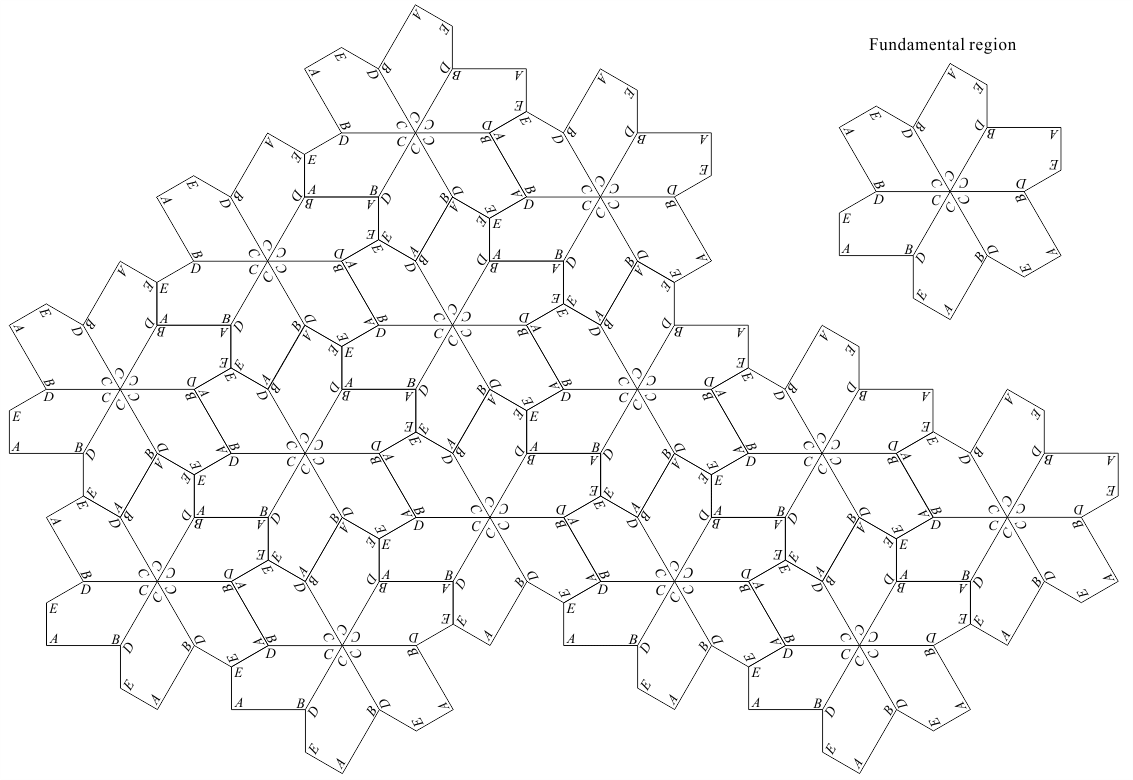} 
  \caption{{\small 
Representative tiling of Type 5. The fundamental region is a unit 
that can generate a periodic tiling by translation only.} 
\label{fig3}
}
\end{figure}

Convex pentagons in Figure~\ref{fig1} can be divided into two equilateral triangles 
\textit{ABD} and \textit{BCD}, and a isosceles triangle \textit{ADE}. The 
area of isosceles triangle \textit{ADE} is equal to 1/3 of the area of 
equilateral triangle \textit{ABD}. That is, the convex pentagons 
in Figure~\ref{fig1} include 7/3 equilateral triangles. As will be described later, 
the convex pentagon in Figure~\ref{fig1} can be considered as a unique convex 
pentagon obtained from a \textbf{t}risected \textbf{h}eptiamond\footnote{ A polyiamond (simply 
iamond) is a plane figure constructed by congruent equilateral triangle 
joined edge-to-edge~\cite{G_and_S_1987, wiki-15, wiki-17}. A polyiamond with seven 
equilateral triangles is called a heptiamond or 7-iamonds. There are 24 unique 
pieces~\cite{G_and_S_1987}. On the other hand, the convex pentagon in Figure~\ref{fig1} is also 
14-polydrafter~\cite{wiki-14, wiki-15}. However, the convex pentagon corresponding to a 
14-polydrafter also exists other than the pentagon in Figure~\ref{fig1}. It belongs 
to Types 1, 5, and 6 (see Figure 8(p) in \cite{Sugimoto_2016}).}. Hereafter, the convex 
pentagon of Figure~\ref{fig1} is referred to as a \textit{TH-pentagon}. The TH-pentagon 
belongs to both Type 1 and Type 5 of the known types for convex pentagonal 
tiles (see Figure~\ref{fig55} in Appendix)~\cite{G_and_S_1987, Sugimoto_2015, wiki-13}. 
Therefore, as shown in Figures~\ref{fig2} and \ref{fig3}, the TH-convex 
pentagon can generate each representative tiling of Type 1 or Type 5, or variations 
of Type 1 tilings (i.e., tilings whose vertices are formed only by the relations 
of $A + D + E = 360^ \circ $ and $B + C = 180^ \circ$).

\renewcommand{\figurename}{{\small Figure.}}

\renewcommand{\figurename}{{\small Figure.}}

\renewcommand{\figurename}{{\small Figure.}}

\begin{figure}[htbp]
 \centering\includegraphics[width=15cm,clip]{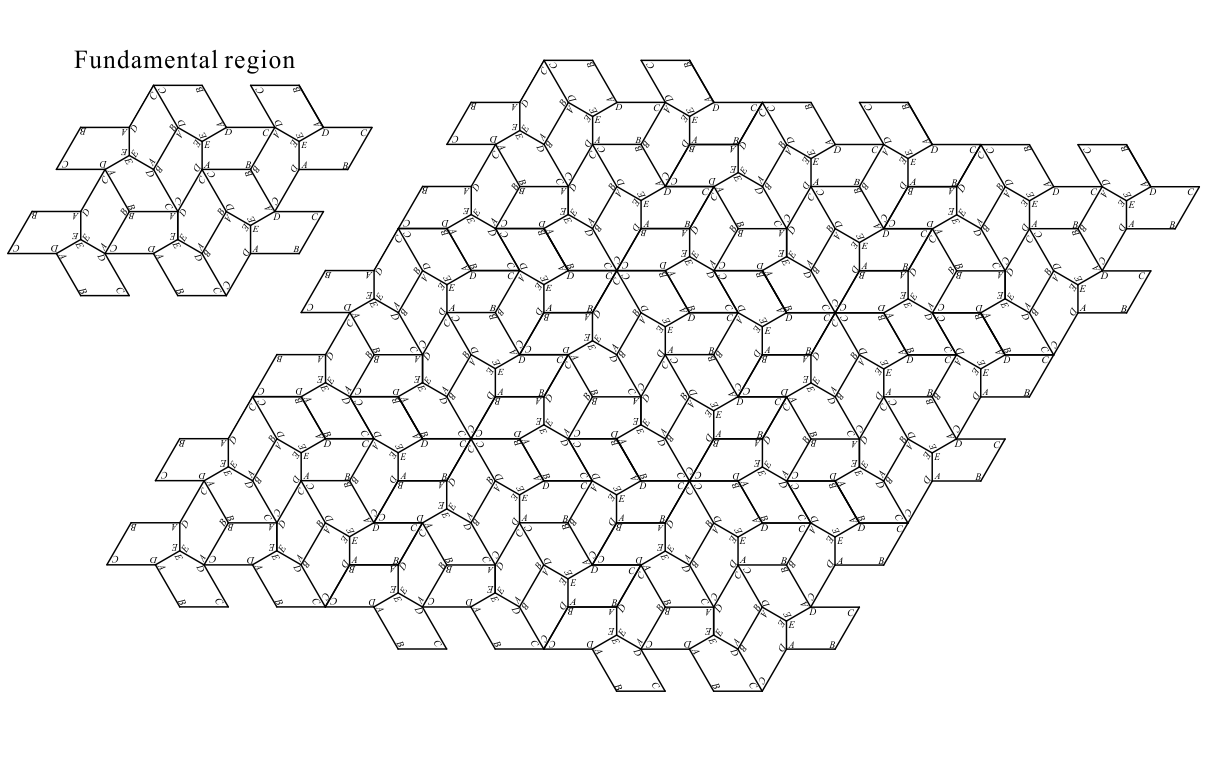} 
  \caption{{\small 
Rice1995-tiling.} 
\label{fig4}
}
\end{figure}

In 1995, Marjorie Rice discovered the interesting tiling in Figure~\ref{fig4} by 
using the TH-pentagon~\cite{Peterson_2010, Scherphuis, Sugimoto_2011, wiki-13, MAA}.
The tiling in Figure~\ref{fig4} decorates the floor of the Marcia P Sward 
lobby in the Mathematical Association of America headquarters. Hereafter, 
the tiling discovered by Marjorie Rice in 1995 is called a \textit{Rice1995-tiling}. 
The fundamental region in Figure~\ref{fig4} is formed by 18 convex pentagons.

\begin{figure}[htbp]
 \centering\includegraphics[width=3cm,clip]{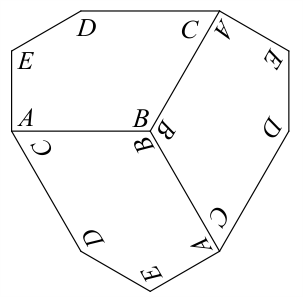} 
  \caption{{\small 
Convex nonagon unit (CN-unit).} 
\label{fig5}
}
\end{figure}

\begin{figure}[htbp]
 \centering\includegraphics[width=14cm,clip]{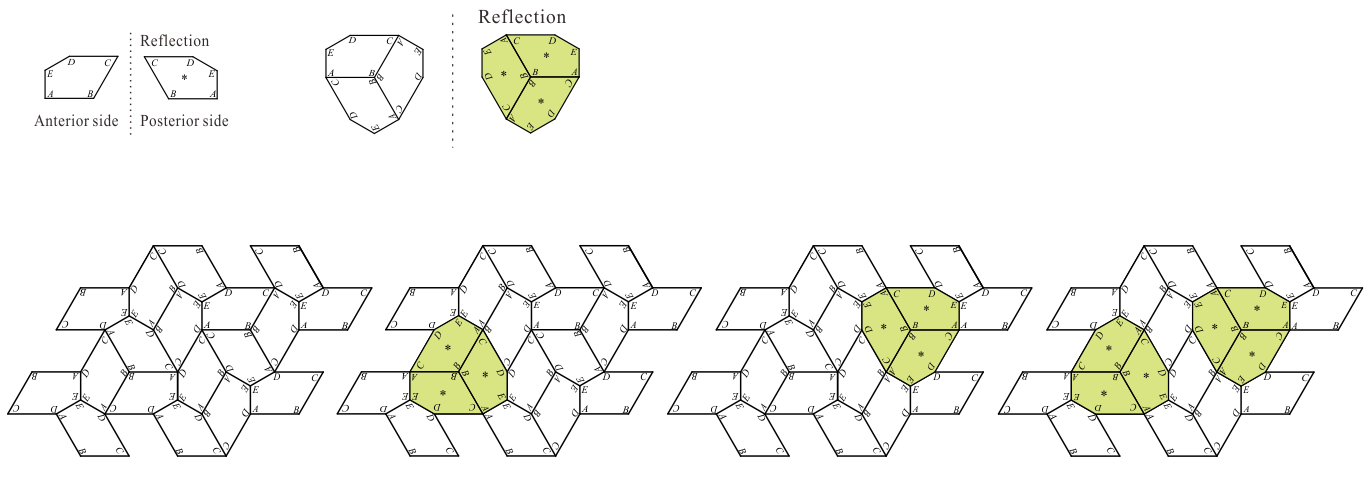} 
  \caption{{\small 
Fundamental regions of a Rice1995-tiling and convex nonagon units.} 
\label{fig6}
}
\end{figure}

\begin{figure}[htbp]
 \centering\includegraphics[width=14cm,clip]{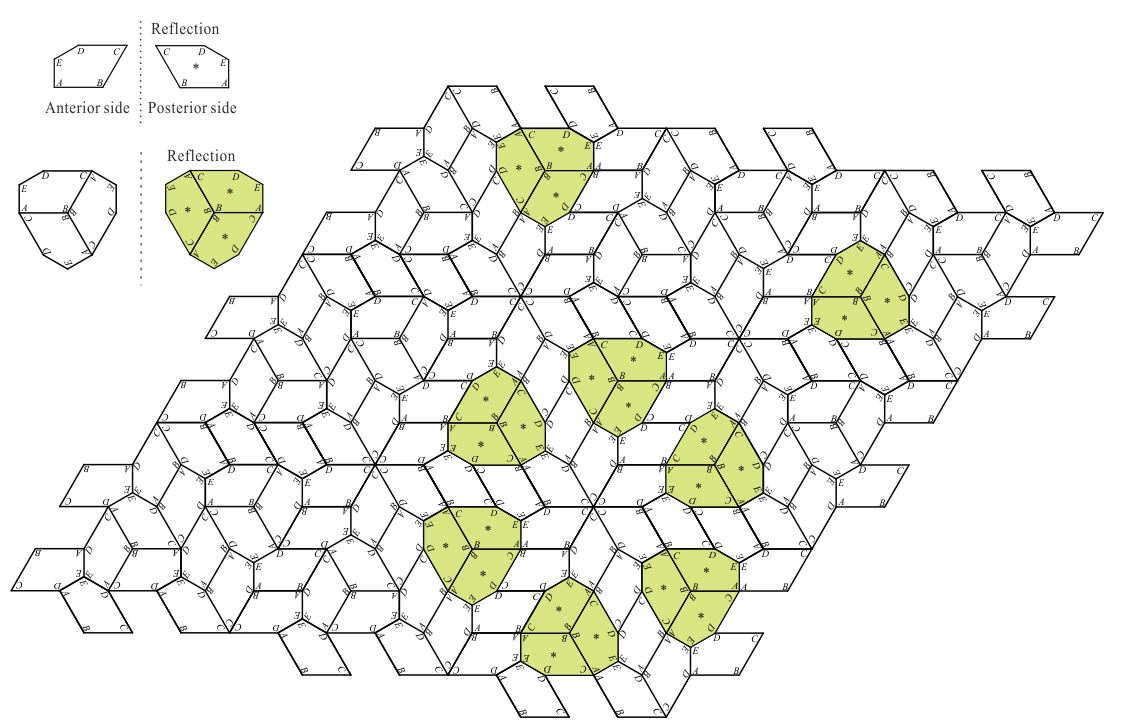} 
  \caption{{\small 
Example of a Rice1995-tiling which contains reflection of convex 
nonagon units.} 
\label{fig7}
}
\end{figure}

As shown in Figure~\ref{fig5}, the authors consider a convex nonagon that is formed 
by three TH-pentagons. The convex nonagon in Figure~\ref{fig5} is called a 
\textit{convex nonagon unit} (CN-unit). 
There are two convex nonagon units in the fundamental region of Figure~\ref{fig4}. 
The outline of a convex nonagon unit is symmetrical (i.e., it is identical 
to the original when reflected about an axis). Therefore, when the convex 
nonagon units in the fundamental region of a Rice1995-tiling are reflected, 
the outline of fundamental region of a Rice1995-tiling is the same (see 
Figure~\ref{fig6}). As shown in Figure~\ref{fig7}, the units in a convex nonagon in 
Rice1995-tiling can be reversed freely~\cite{Scherphuis}. Thus, from the property 
for convex nonagon units, the Rice1995-tiling has the property of generating 
nonperiodic tilings.

Although the TH-pentagon is known to be able to generate tilings of Type 1, 
Type 5, or Rice1995, the authors discovered that it is possible to generate 
tilings other than these. The newly discovered tilings (tessellations) are 
presented in this manuscript. Note that many proofs are omitted.

\renewcommand{\figurename}{{\small Figure.}}

\renewcommand{\figurename}{{\small Figure.}}

\renewcommand{\figurename}{{\small Figure.}}

\renewcommand{\figurename}{{\small Figure.}}

\section{How to concentrate vertices}
\label{section2}

If the TH-pentagon generates tilings, it must have the concentration 
relation that the sum of some vertices (angles) is equal either to 
$180^ \circ$ or to $360^ \circ$. First, the relations of $180^ \circ$ are 
combinations of the following.

\medskip\noindent
$
2A = 180^ \circ ,
\quad
B + C = 180^ \circ ,
\quad
B + E = 180^ \circ ,
\quad
C + E = 180^ \circ .
$
\medskip

\noindent
Next, the concentration relations of $360^ \circ$ are combinations of the 
following.

\medskip\noindent
$
A + B + D = 360^ \circ ,
\quad
A + D + E = 360^ \circ ,
\quad
3B = 360^ \circ ,
\quad
2B + E = 360^ \circ ,
\quad
B + 2E = 360^ \circ ,
\quad
C + 2D = 360^ \circ ,
\quad
B + 2E = 360^ \circ ,
\quad
3E = 360^ \circ ,
\quad
4A = 360^ \circ ,
\quad
2A + B + C = 360^ \circ ,
\quad 
2A + C + E = 360^ \circ ,
\quad
A + 2C + D = 360^ \circ ,
\quad
2B + 2C = 360^ \circ ,
\quad
B + 2C + E = 360^ \circ ,
\quad
2C + 2E = 360^ \circ ,
\quad
2A + 3C = 360^ \circ ,
\quad
B + 4C = 360^ \circ ,
\quad
4C + E = 360^ \circ ,
\quad
6C = 360^ \circ .
$

\begin{figure}[htbp]
 \centering\includegraphics[width=11cm,clip]{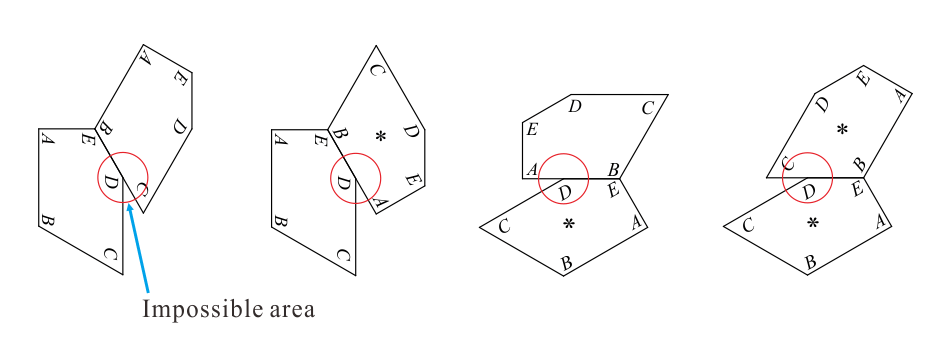} 
  \caption{{\small 
Example of impossible areas where TH-pentagon cannot be aligned 
without gaps.} 
\label{fig8}
}
\end{figure}

\medskip\noindent
However, if $B + E = 180^ \circ $, impossible areas are created where a 
TH-pentagon cannot be aligned without gaps always exist (see Figure~\ref{fig8}). 
The concentration relations that cannot be used in tilings clearly are excluded. 
As a result, candidates of concentration relations that will be used in 
tilings are as follows.

\medskip\noindent
$
2A = 180^ \circ ,
\quad
B + C = 180^ \circ ,
\quad\\
A + B + D = 360^ \circ ,
\quad
A + D + E = 360^ \circ ,
\quad
3B = 360^ \circ ,
\quad
C + 2D = 360^ \circ ,
\quad
3E = 360^ \circ ,
\quad
4A = 360^ \circ ,
\quad
2A + B + C = 360^ \circ ,
\quad
2A + C + E = 360^ \circ ,
\quad
A + 2C + D = 360^ \circ ,
\quad
2B + 2C = 360^ \circ ,
\quad
2A + 3C = 360^ \circ ,
\quad
B + 4C = 360^ \circ ,
\quad
6C = 360^ \circ .
$
\medskip

\begin{figure}[htbp]
 \centering\includegraphics[width=13cm,clip]{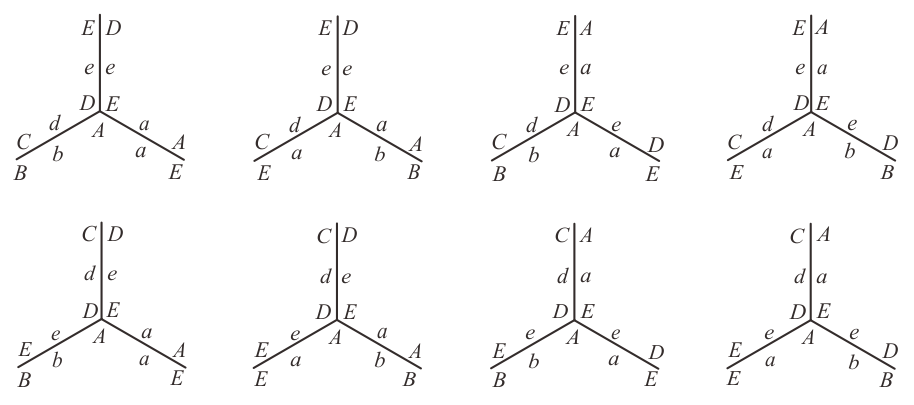} 
  \caption{{\small 
Eight possible arrangements of $A + D + E = 360^ \circ$.} 
\label{fig9}
}
\end{figure}

\begin{figure}[htbp]
 \centering\includegraphics[width=14cm,clip]{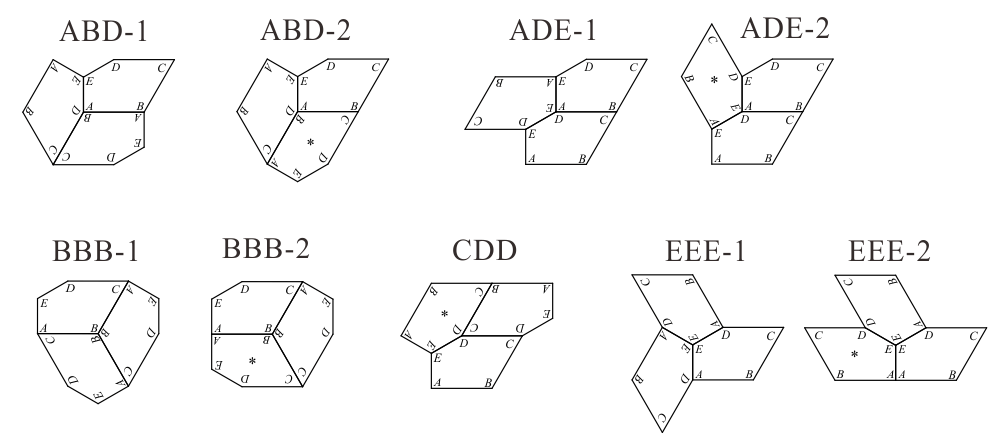} 
  \caption{{\small 
Unique patterns of concentration relations with three vertices.} 
\label{fig10}
}
\end{figure}

\begin{figure}[htbp]
 \centering\includegraphics[width=15cm,clip]{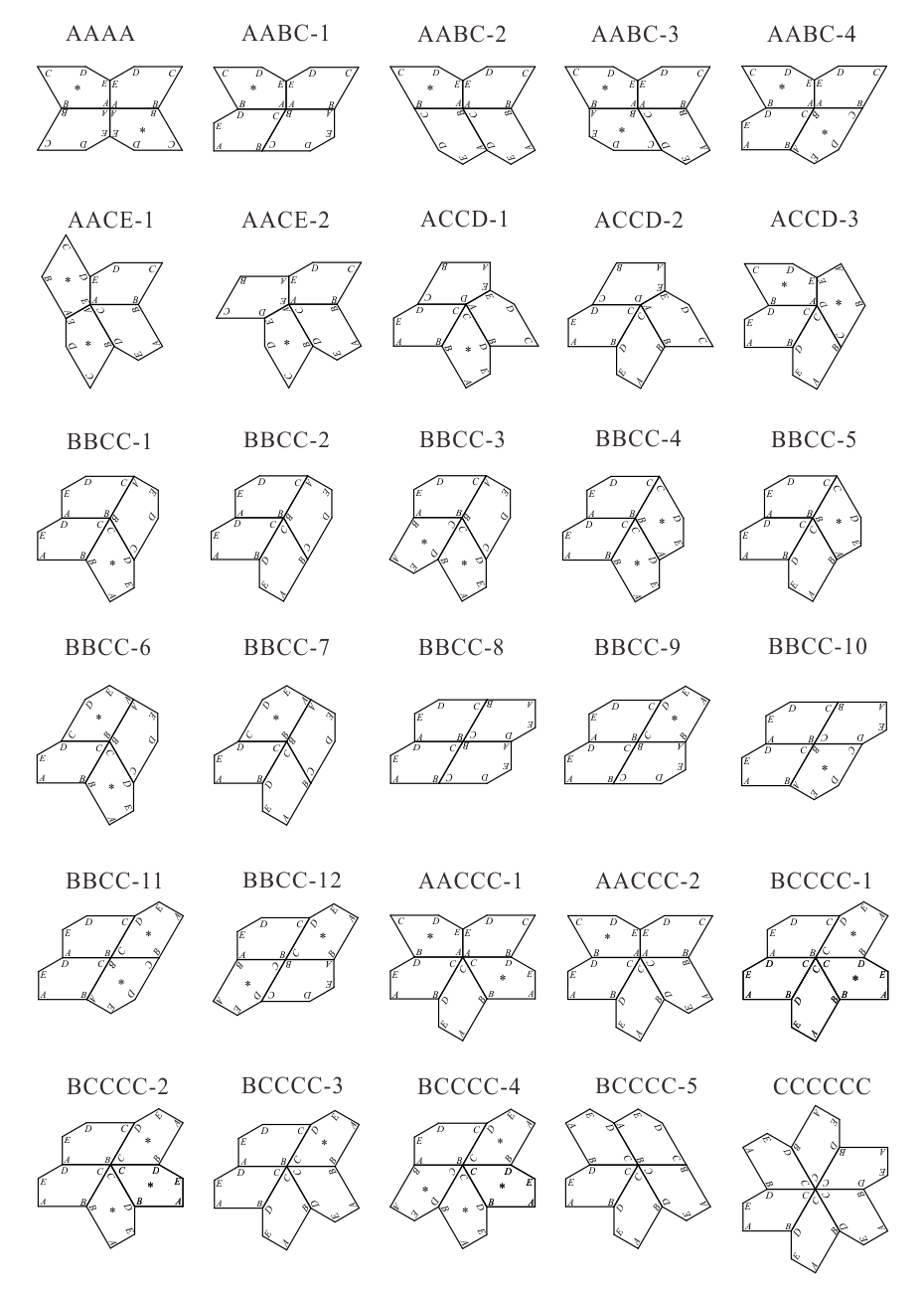} 
  \caption{{\small 
Unique patterns of concentration relations with four or more 
vertices.} 
\label{fig11}
}
\end{figure}

Next, consider the arrangements of TH-pentagons in these candidates of 
concentration relations. For example, for concentration relations with three 
vertices like $A + D + E = 360^ \circ$, there are eight possible arrangements 
to assemble the three pentagons since each of three pentagons can be 
reversed (see Figure~\ref{fig9}). Except for arrangements that clearly cannot 
be used in tilings, the arrangements of concentration relations with three 
vertices are nine patterns in Figure~\ref{fig10}, not distinguishing reflections and 
rotations (i.e., the number of unique patterns in this case is nine). The 
other cases are also obtained similarly (see Figure~\ref{fig11}). For example, for 
$2A + C + E = 360^ \circ$, two patterns AACE-1 and AACE-2 as shown in 
Figure~\ref{fig11} are obtained.

Here, consider the concentration $A + D + E = 360^ \circ $, $3E = 360^ \circ $, 
and $2A + C + E = 360^ \circ $ including the vertex $E$. In the case of 
$2A + C + E = 360^ \circ $, there are two patterns AACE-1 and AACE-2 
in Figure~\ref{fig11}, but it is confirmed that an impossible area appears by 
accumulating the step by step possible arrangements of TH-pentagons. 
Therefore, if the TH-pentagon generates a tiling, the vertex $E$ belongs to 
$A + D + E = 360^ \circ$ or $3E = 360^ \circ$. Furthermore, the authors 
confirmed that $A + D + E = 360^ \circ$ and $3E = 360^ \circ$ cannot 
coexist within the tiling. 

\begin{figure}[htbp]
 \centering\includegraphics[width=14cm,clip]{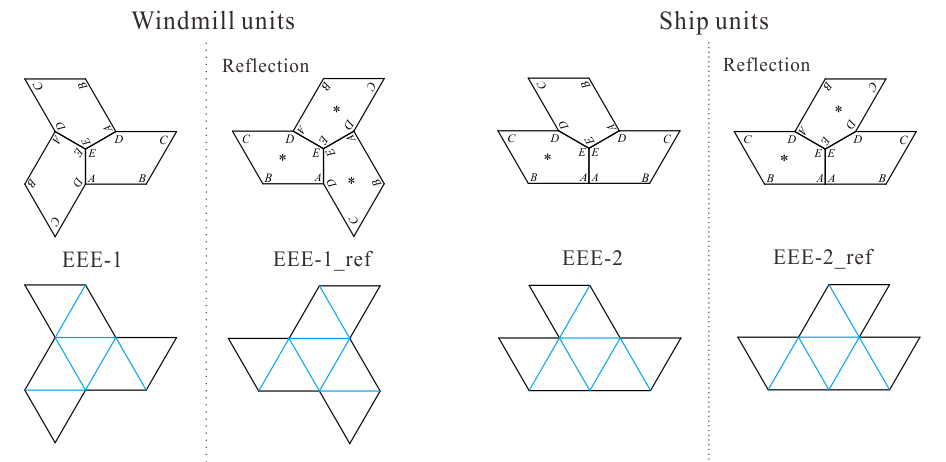} 
  \caption{{\small 
Windmill units and ship units.} 
\label{fig12}
}
\end{figure}

Therefore, the TH-pentagon can generate tilings with $A + D + E = 360^ \circ$ 
and tilings with $3E = 360^ \circ $. The tilings with $A + D + E = 360^ \circ$ 
are always variations of Type 1 tilings as shown in Figure~\ref{fig2}. The 
tilings with $3E = 360^ \circ $ contain the patterns EEE-1 or EEE-2 as shown 
Figure~\ref{fig10}. Hereafter, EEE-1 and EEE-2 are referred to as a 
\textit{windmill unit} and a \textit{ship unit}, respectively. 
As shown in Figure~\ref{fig12}, the windmill unit and the ship unit can 
be considered as heptiamonds. That is, the TH-pentagon is a convex pentagon 
that can be obtained by dividing the two types of heptiamond of the windmill 
unit and the ship unit into three equal parts. The convex pentagon that can 
be obtained by trisecting 24 types of heptiamond does not exist other than 
the TH-pentagon. Therefore, tilings with $3E = 360^ \circ$ by TH-pentagons 
are equivalent to tilings of these two types of heptiamonds.

The authors discovered new convex pentagon tilings using windmill units and 
ship units; they are presented in the following sections.

\renewcommand{\figurename}{{\small Figure.}}

\renewcommand{\figurename}{{\small Figure.}}

\renewcommand{\figurename}{{\small Figure.}}

\renewcommand{\figurename}{{\small Figure.}}

\renewcommand{\figurename}{{\small Figure.}}

\section{Tilings by only the windmill units}
\label{section3}

\begin{figure}[htbp]
 \centering\includegraphics[width=8cm,clip]{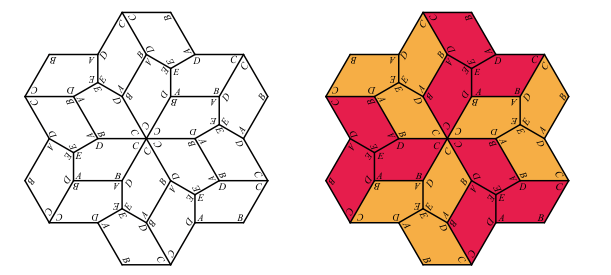} 
  \caption{{\small 
Hexagonal flowers L1 unit.} 
\label{fig13}
}
\end{figure}

In this section, tilings by only the windmill units are introduced. 
As shown in Figure~\ref{fig13}, the authors consider a shape (18 sided 
polygon) formed by six windmill units. Hereafter, the shape is referred 
to as a \textit{hexagonal flowers L1 unit} (HFL1-unit). The hexagonal 
flowers L1 unit has an edge in common with edge $b$, $c$, or $d$ of a 
TH-pentagon.

\begin{figure}[htbp]
 \centering\includegraphics[width=14.5cm,clip]{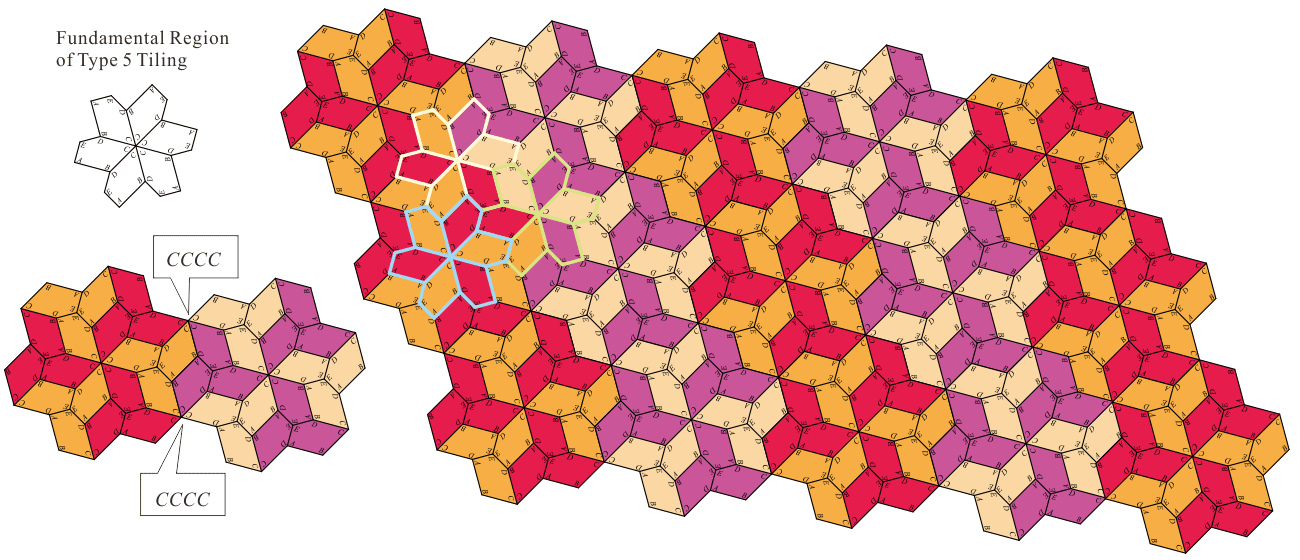} 
  \caption{{\small 
Relation of representative Type 5 tiling and hexagonal flower L1 
units. Note that, in order to make a contiguity method intelligible, the 
hexagonal flowers L1 units of different colors are used.} 
\label{fig14}
}
\end{figure}

\begin{figure}[htbp]
 \centering\includegraphics[width=14.5cm,clip]{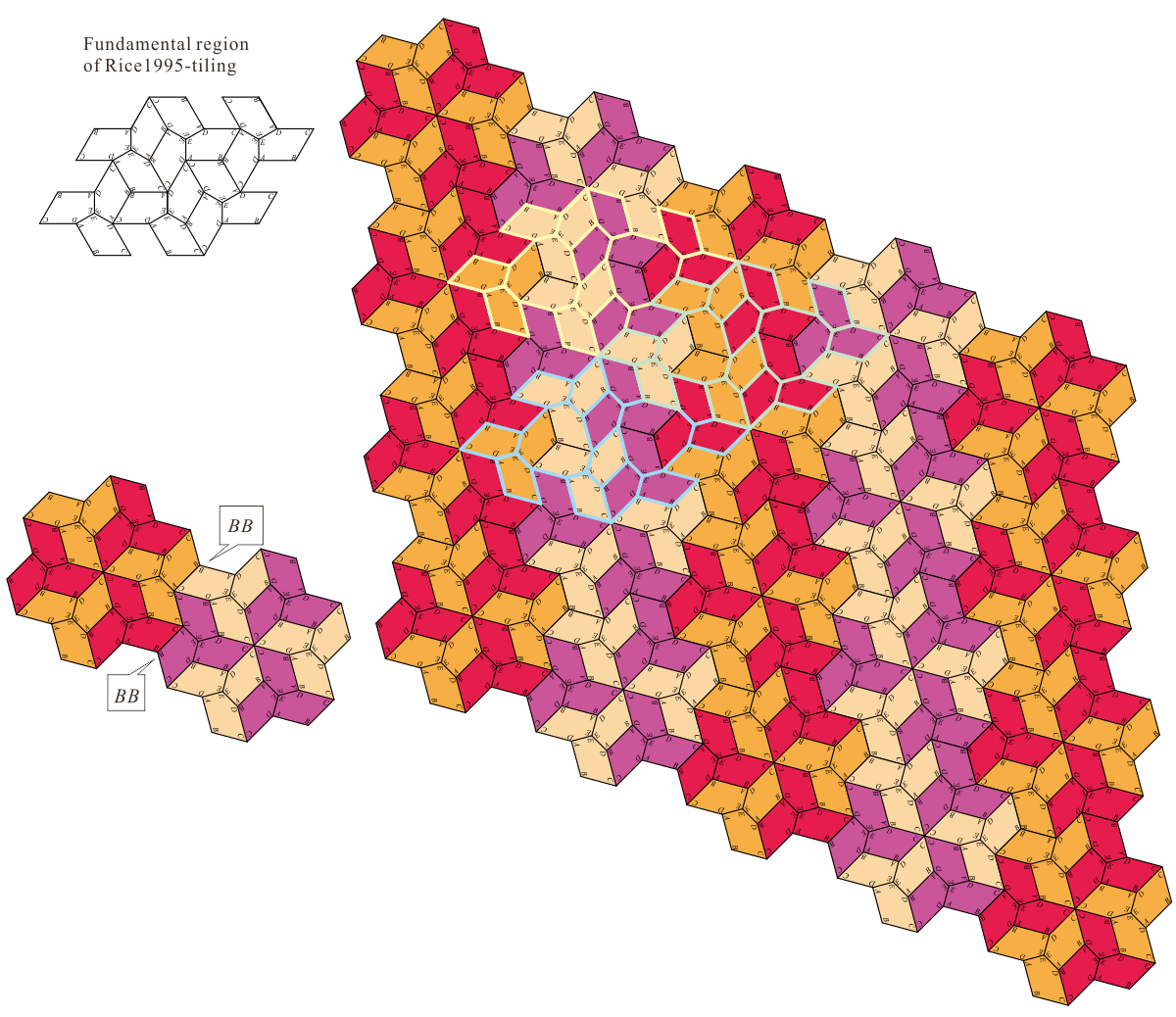} 
  \caption{{\small 
Relation of Rice1995-tiling and hexagonal flowers L1 units. Note 
that, in order to make a contiguity method intelligible, the hexagonal 
flowers L1 units of different colors are used.} 
\label{fig15}
}
\end{figure}

The representative Type 5 tiling (Figure~\ref{fig3}) and the 
Rice1995-tiling (Figure~\ref{fig4}) can be generated using the hexagonal 
flowers L1 units. As shown in Figures~\ref{fig14} and \ref{fig15}, the 
representative Type 5 tiling and the Rice1995-tiling that does not 
contain reflected convex nonagon units are generated by the 
difference in the contiguity method of hexagonal flowers L1 units.

\renewcommand{\figurename}{{\small Figure.}}

\renewcommand{\figurename}{{\small Figure.}}

\renewcommand{\figurename}{{\small Figure.}}

\begin{figure}[htbp]
 \centering\includegraphics[width=13cm,clip]{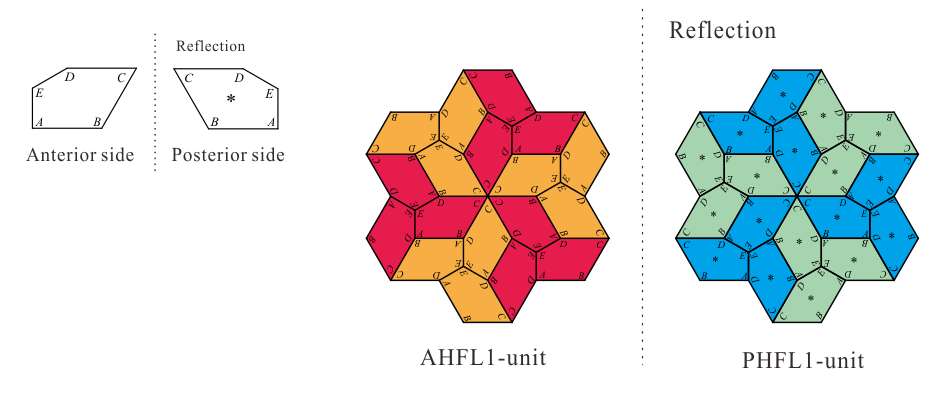} 
  \caption{{\small 
AHFL1-unit (anterior hexagonal flowers L1 unit) and PHFL1-unit 
(posterior hexagonal flowers L1 unit).} 
\label{fig16}
}
\end{figure}

The hexagonal flowers L1 unit and the reflection of the hexagonal flowers L1 
unit have identical outlines (see Figure~\ref{fig16}). The eighteen convex pentagons 
in the hexagonal flowers L1 units in Figure~\ref{fig13} are anterior since the 
pentagon in Figure~\ref{fig1} is considered to be forward facing. Hereafter, the 
hexagonal flowers L1 unit formed by 18 anterior TH-pentagons is referred to 
as an \textit{AHFL1-unit}, and the hexagonal flowers L1 unit that is formed by 
18 posterior TH-pentagons is referred to as a \textit{PHFL1-unit}. In addition, 
the convex nonagon unit (see Figure~\ref{fig5}) formed by three anterior 
TH-pentagons is referred to as an \textit{ACN-unit}, and the convex nonagon 
unit formed by three posterior TH-pentagons is referred to as a \textit{PCN-unit}.

\begin{figure}[htbp]
 \centering\includegraphics[width=15.5cm,clip]{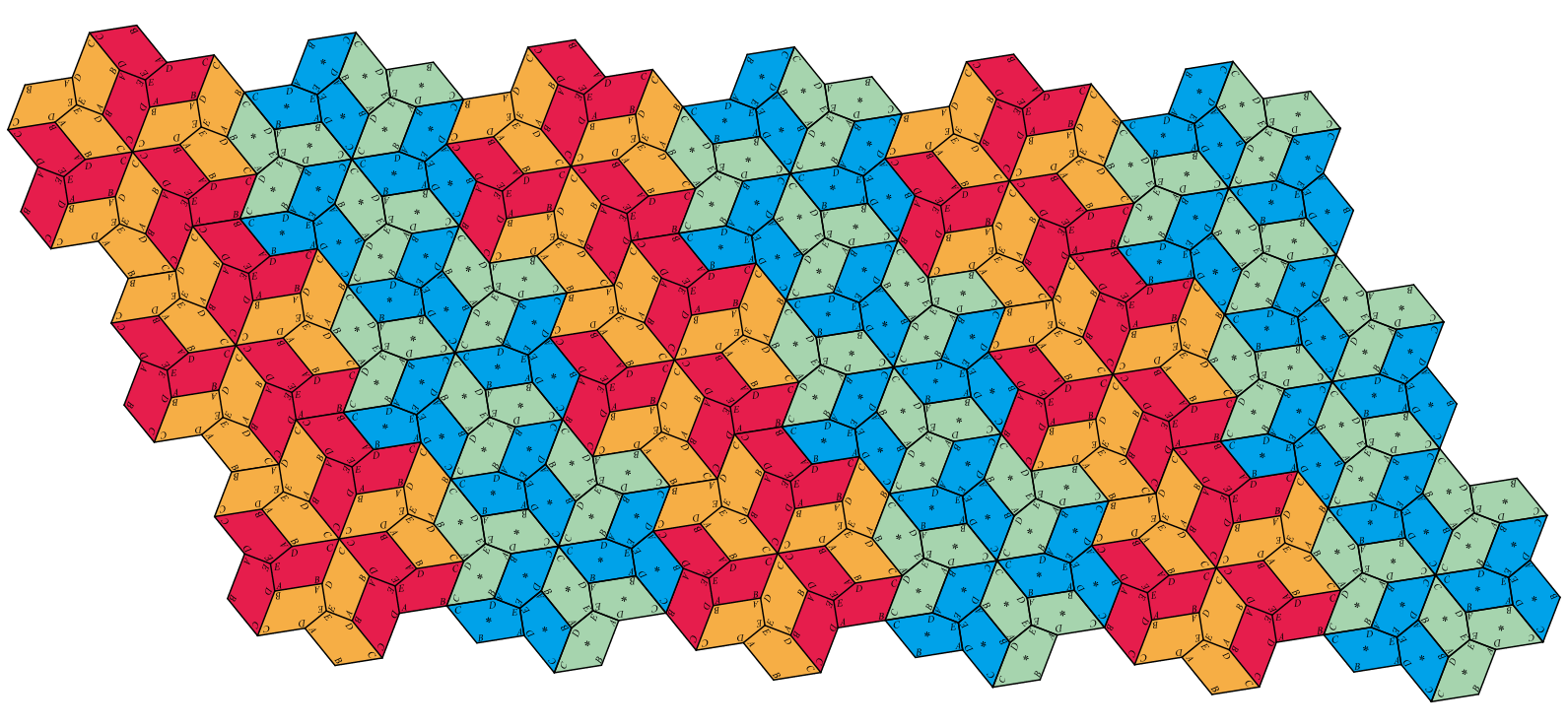} 
  \caption{{\small 
Tiling by a pair of an AHFL1-unit and a PHFL1-unit.} 
\label{fig17}
}
\end{figure}

\begin{figure}[htbp]
 \centering\includegraphics[width=8cm,clip]{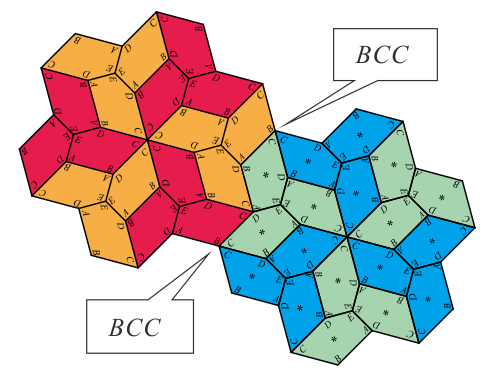} 
  \caption{{\small 
Contiguity method of an AHFL1-unit and a PHFL1-unit.} 
\label{fig18}
}
\end{figure}

Figure~\ref{fig17} is a tiling consisting of a pair of AHFL1-unit and PHFL1-unit. 
(The tilings in Figure~\ref{fig17} have 5-valent vertices ``$B + 4C = 360^ \circ$.'') 
The contiguity method of AHFL1-units and PHFL1-units is only a pattern 
as shown in Figure~\ref{fig18}. 

\begin{figure}[htbp]
 \centering\includegraphics[width=15.5cm,clip]{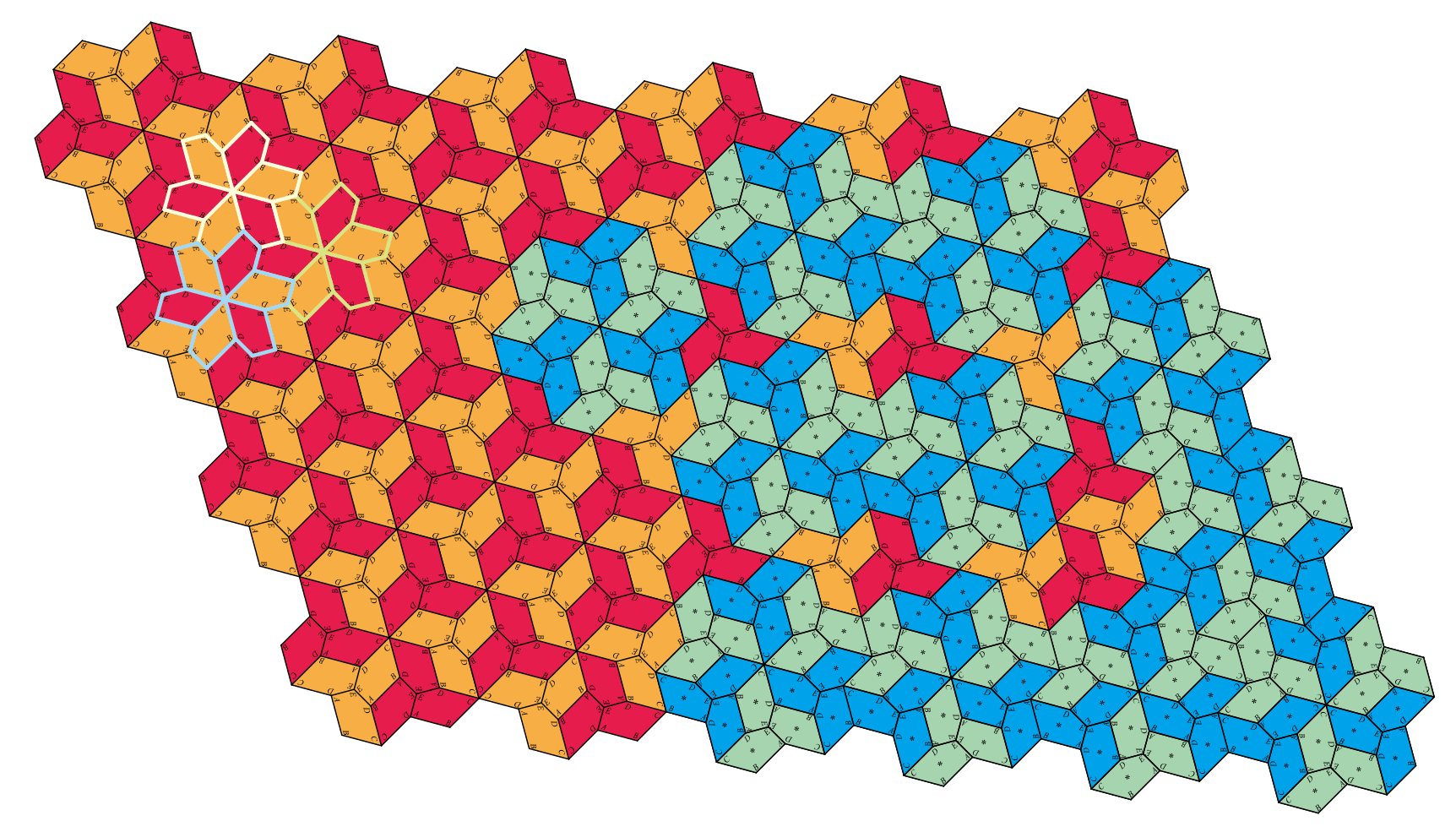} 
  \caption{{\small 
Tiling that incorporates PHFL1-units as the basis a Type 5 tiling 
that is formed by AHFL1-units.} 
\label{fig19}
}
\end{figure}

\begin{figure}[htbp]
 \centering\includegraphics[width=13.5cm,clip]{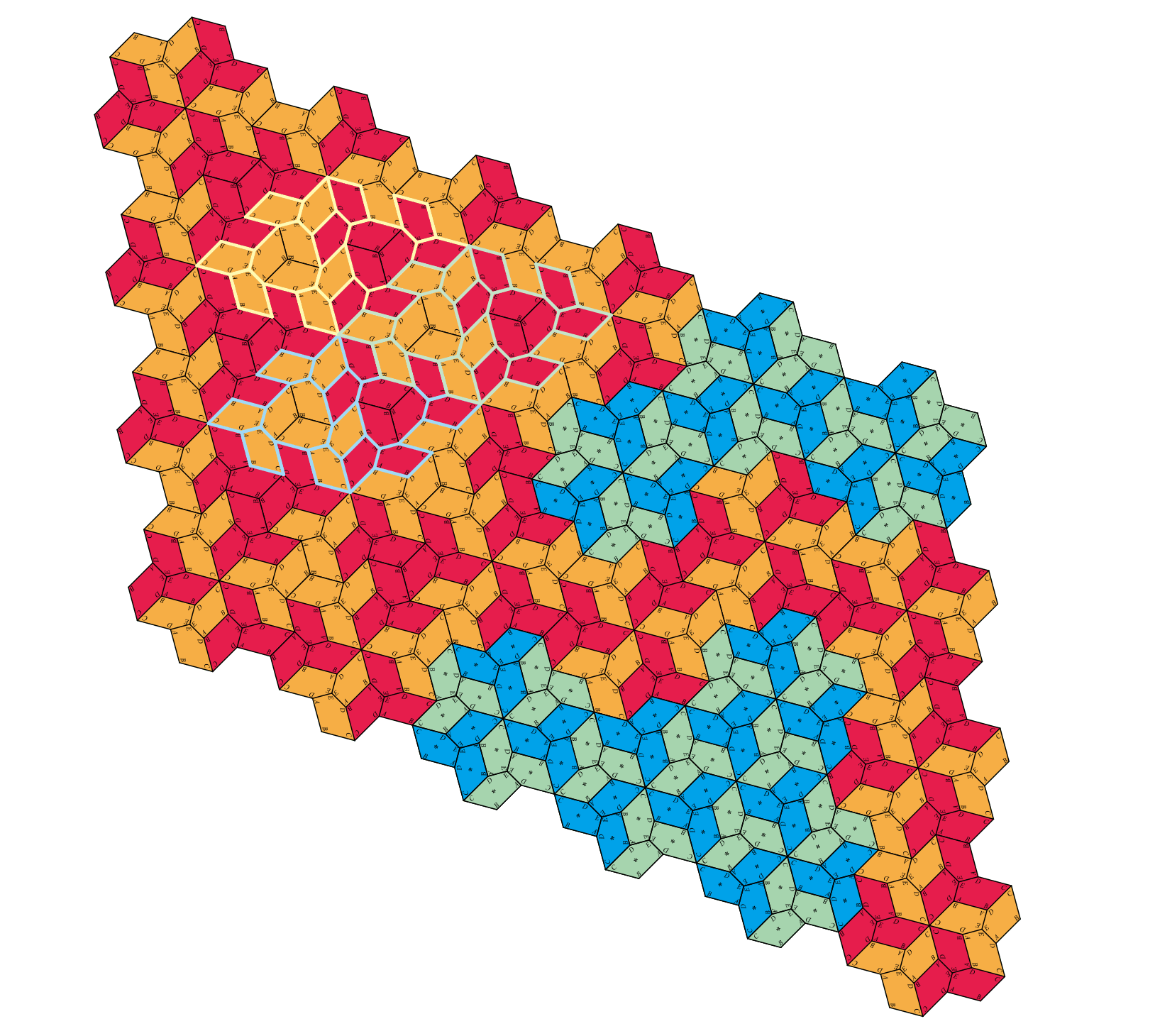} 
  \caption{{\small 
Tiling that incorporates PHFL1-units as the basis of an 
Rice1995-tiling that is formed by AHFL1-units.} 
\label{fig20}
}
\end{figure}

Since the AHFL1-unit and the PHFL1-unit have identical outlines, a hexagonal 
flowers L1 unit in the Type 5 tiling and the Rice1995-tiling can be changed 
to a reversed hexagonal flowers L1 unit freely. Figure~\ref{fig19} is a tiling that 
incorporates PHFL1-units as the basis of a Type 5 tiling that is formed by 
AHFL1-units. Figure~\ref{fig20} is a tiling that incorporates PHFL1-units as the 
basis of a Rice1995-tiling that is formed by AHFL1-units.

\begin{figure}[htbp]
 \centering\includegraphics[width=13.5cm,clip]{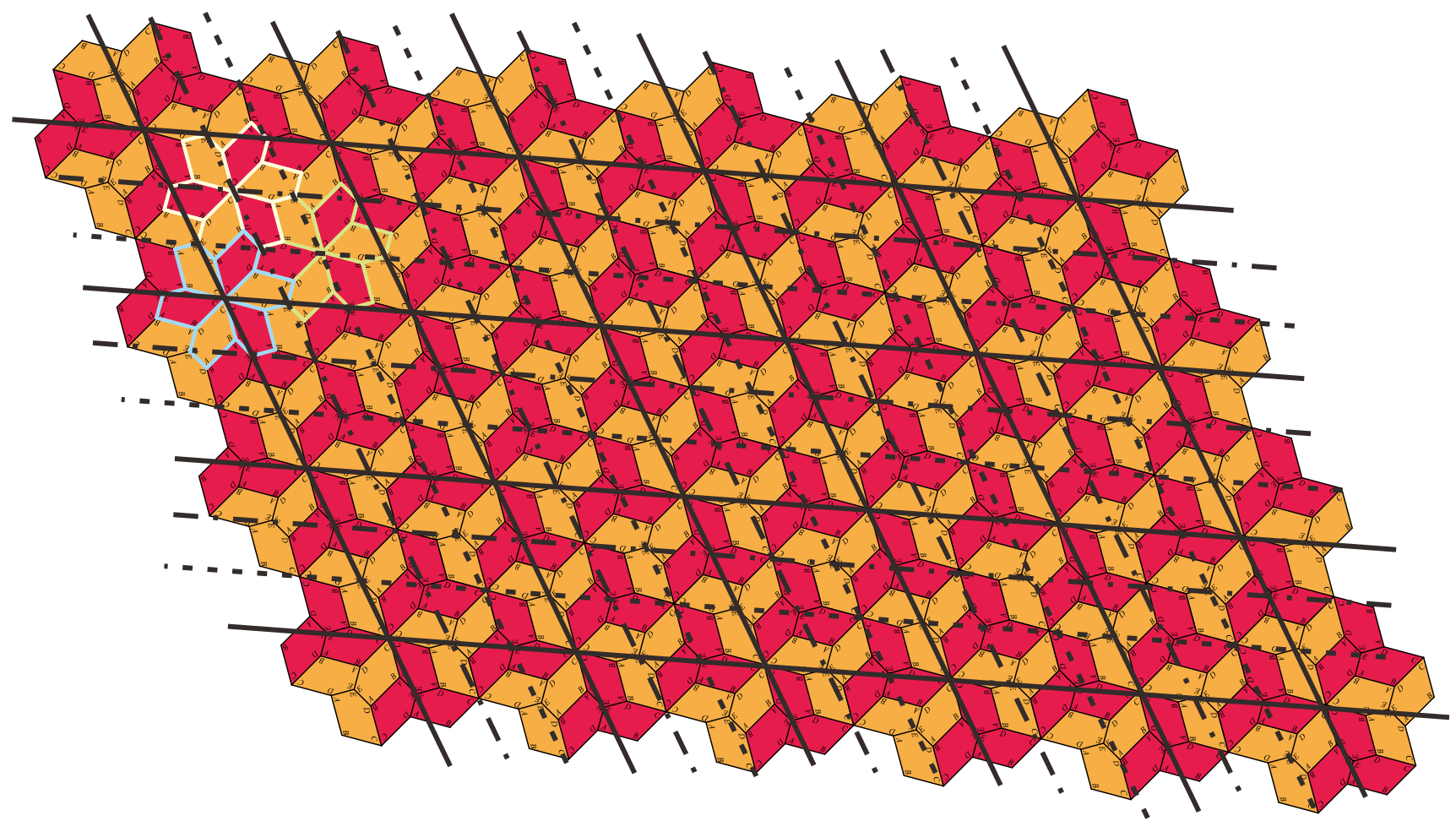} 
  \caption{{\small 
Position relation of points concentrating six vertices $C$ in the 
Type 5 tiling.} 
\label{fig21}
}
\end{figure}

As shown in Figure~\ref{fig21}, points concentrating six vertices $C$ in the 
Type 5 tiling are joined with lines of three kinds (solid line, dashed line, and 
one-dot chain line). The points concentrating six vertices $C$ exist on the 
intersections of the same lines. From Figure~\ref{fig21}, it can be seen that 
there exist three types of points concentrating six vertices $C$.

\renewcommand{\figurename}{{\small Figure.}}

\renewcommand{\figurename}{{\small Figure.}}

\renewcommand{\figurename}{{\small Figure.}}

\renewcommand{\figurename}{{\small Figure.}}

\renewcommand{\figurename}{{\small Figure.}}

\renewcommand{\figurename}{{\small Figure.}}

\renewcommand{\figurename}{{\small Figure.}}
\begin{figure}[htbp]
 \centering\includegraphics[width=13.5cm,clip]{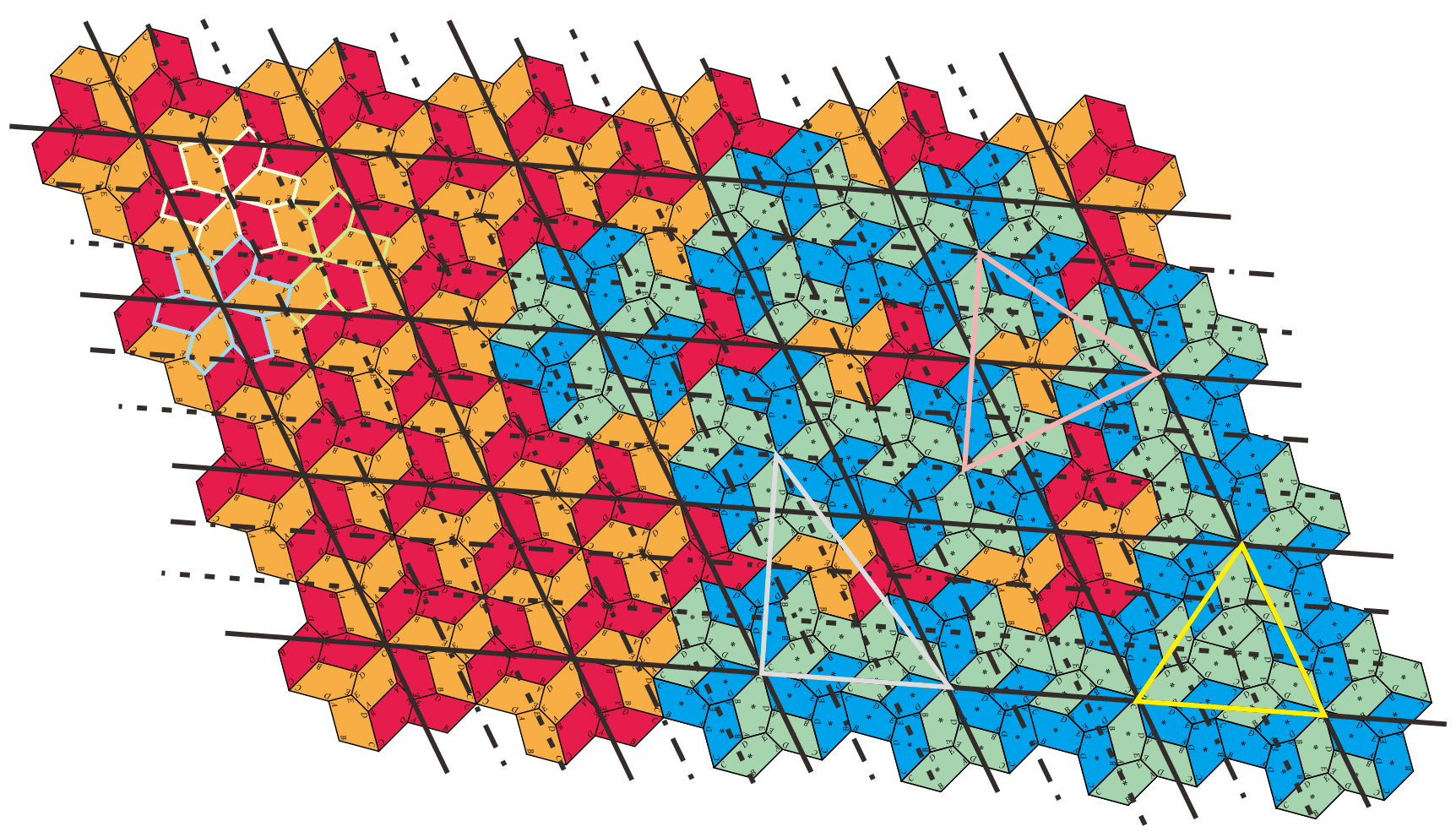} 
  \caption{{\small 
Position relation of points concentrating six vertices $C$ in the 
tiling of Figure~\ref{fig19}.} 
\label{fig22}
}
\end{figure}

Figure~\ref{fig22} transposes the tiling of Figure~\ref{fig19} to the tiling of 
Figure~\ref{fig21}. From Figure~\ref{fig22} it is apparent that three possible 
relations exist between centers of three adjoining PHFL1-units, as follows:

\medskip\noindent
\textbf{Case 1}: The centers of three adjoining PHFL1-units are on the 
intersections of the lines of one kind (see the yellow triangle in Figure~\ref{fig22}).

\medskip\noindent
\textbf{Case 2}: The centers of three adjoining PHFL1-units are on the 
intersections of the lines of two kinds (see the gray triangle in Figure~\ref{fig22}).

\medskip\noindent
\textbf{Case 3}: The centers of three adjoining PHFL1-units are on the 
intersections of the lines of three kinds (see the pink triangle in Figure~\ref{fig22}).
\medskip

\begin{figure}[htbp]
 \centering\includegraphics[width=15cm,clip]{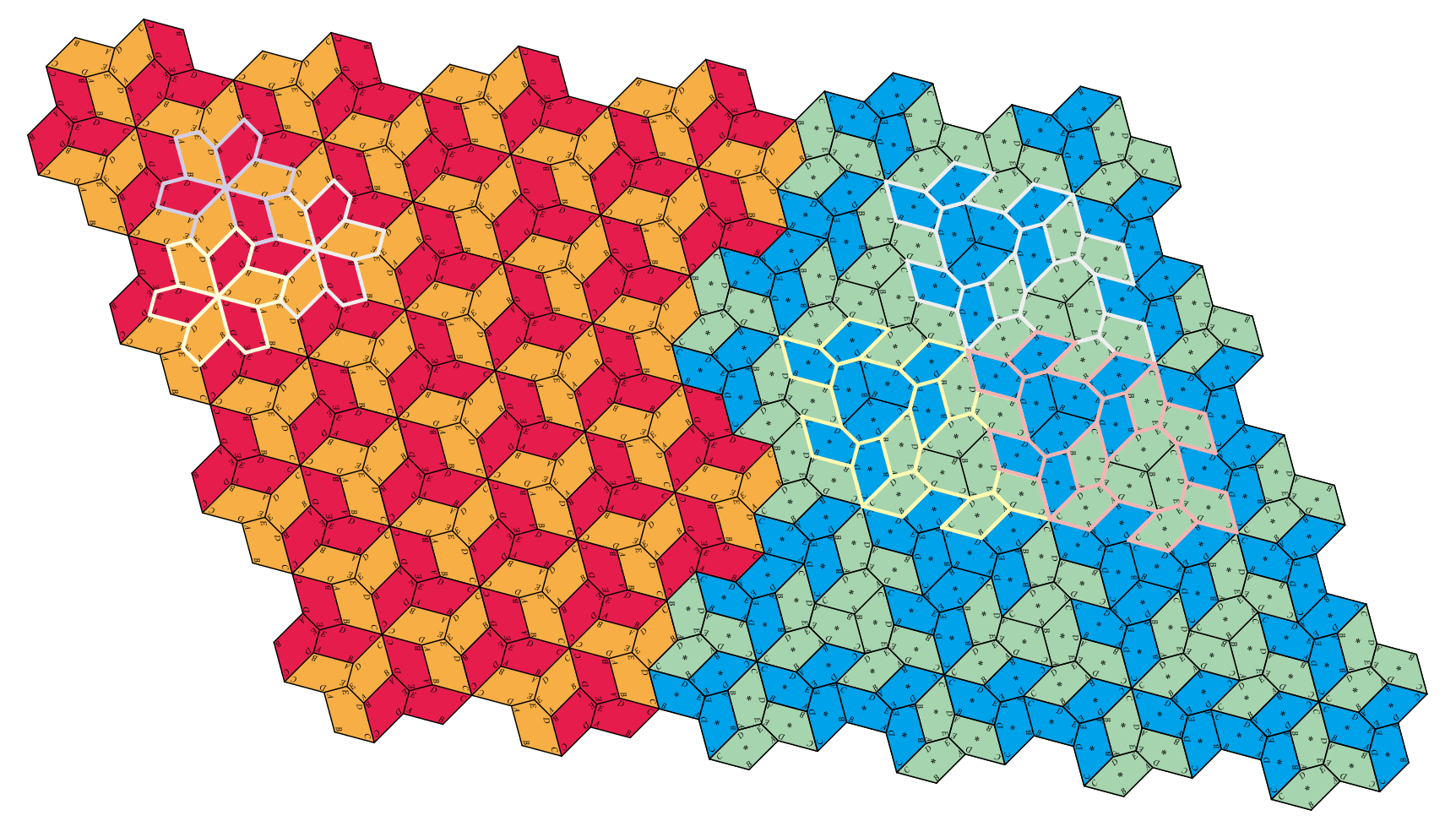} 
  \caption{{\small 
Example of tiling according to the properties of Case 1.} 
\label{fig23}
}
\end{figure}

If a Type 5 tiling with AHFL1-units is filled by PHFL1-units according to 
the properties of Case 1, the Type 5 tiling becomes an Rice1995-tiling with 
PHFL1-units (see Figure~\ref{fig23}). 

\begin{figure}[htbp]
 \centering\includegraphics[width=12cm,clip]{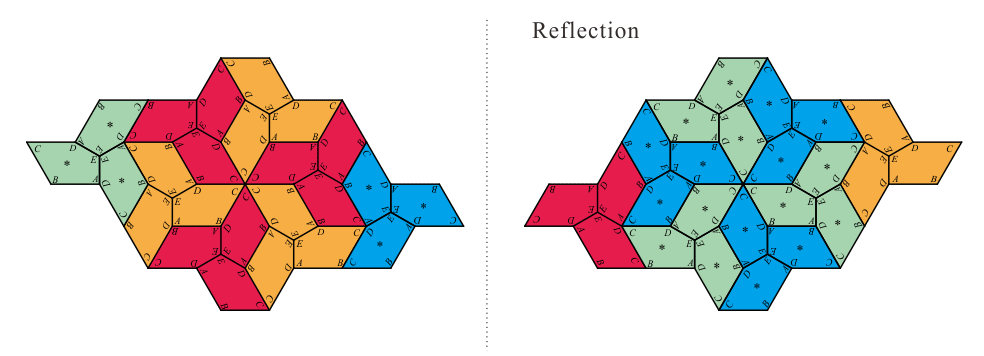} 
  \caption{{\small 
Units that are formed by a hexagonal flowers L1 unit and two 
reverse windmill units.} 
\label{fig24}
}
\end{figure}

\begin{figure}[htbp]
 \centering\includegraphics[width=14cm,clip]{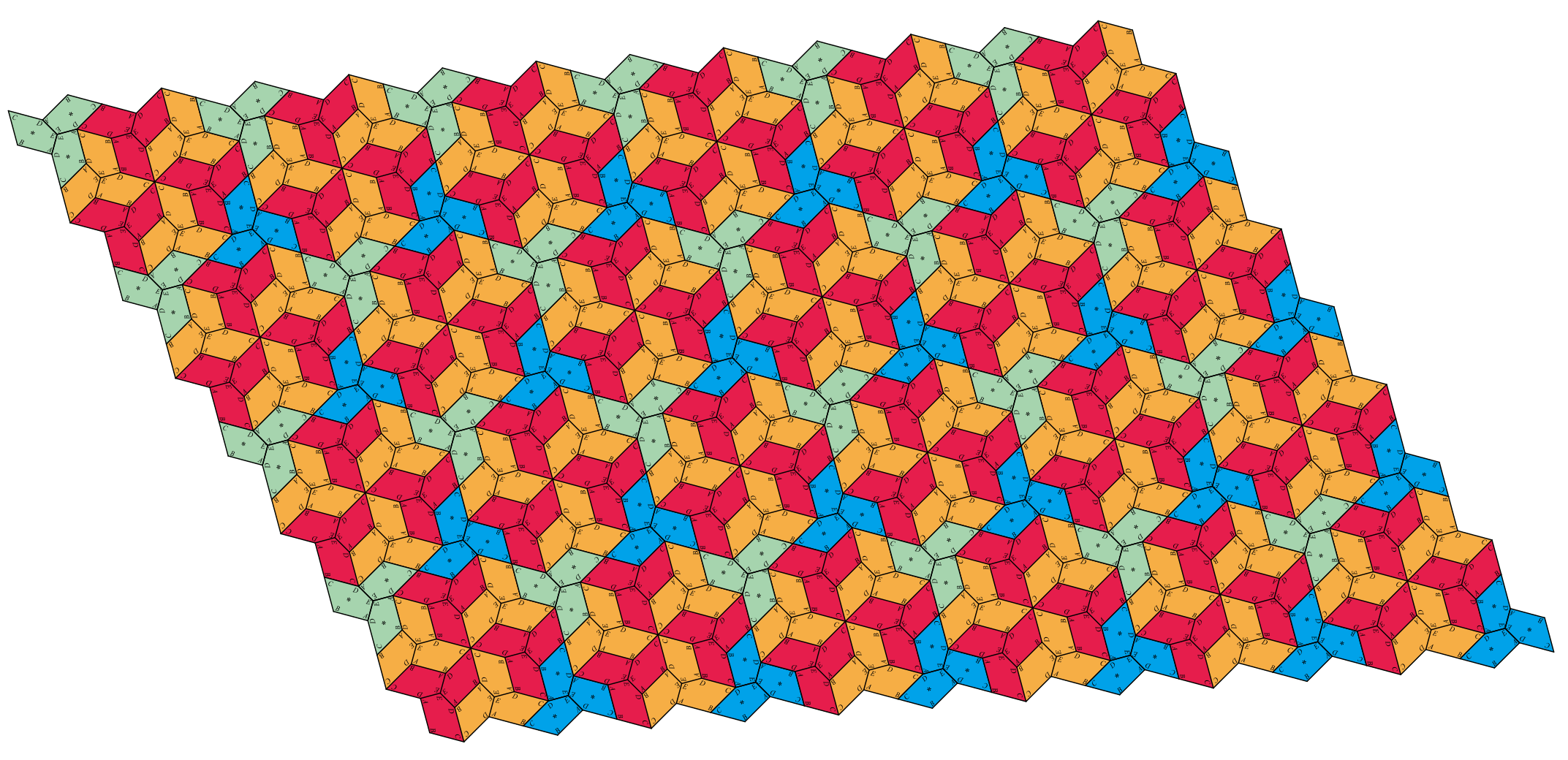} 
  \caption{{\small 
Example of a tiling according to the properties of Case 3.} 
\label{fig25}
}
\end{figure}

\begin{figure}[htbp]
 \centering\includegraphics[width=15cm,clip]{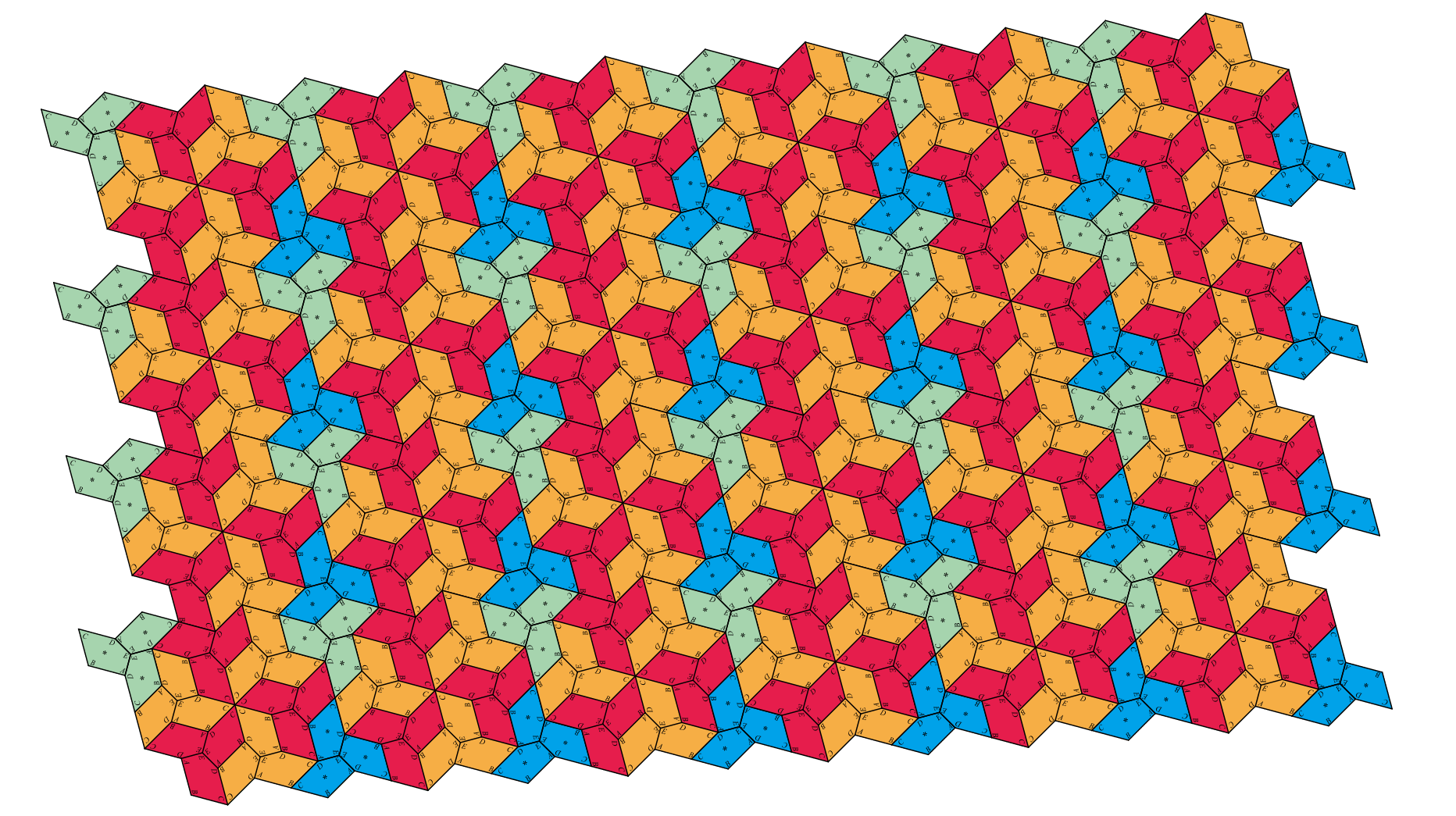} 
  \caption{{\small 
Example of a tiling according to the properties of Case 2.} 
\label{fig26}
}
\end{figure}

\begin{figure}[htbp]
 \centering\includegraphics[width=13.5cm,clip]{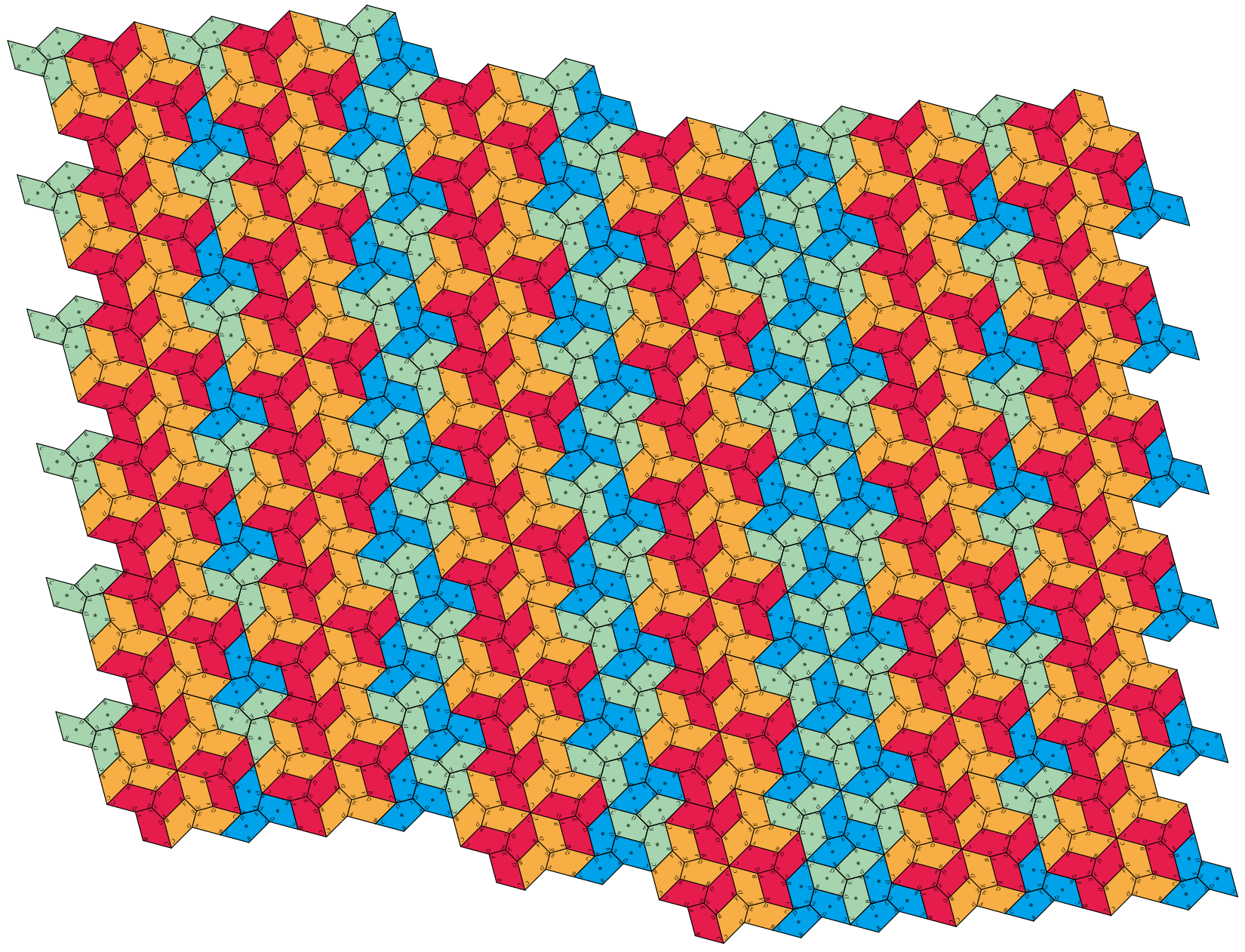} 
  \caption{{\small 
Example of a tiling that is formed by applying the properties of 
Case 2.} 
\label{fig27}
}
\end{figure}

\begin{figure}[htbp]
 \centering\includegraphics[width=13.5cm,clip]{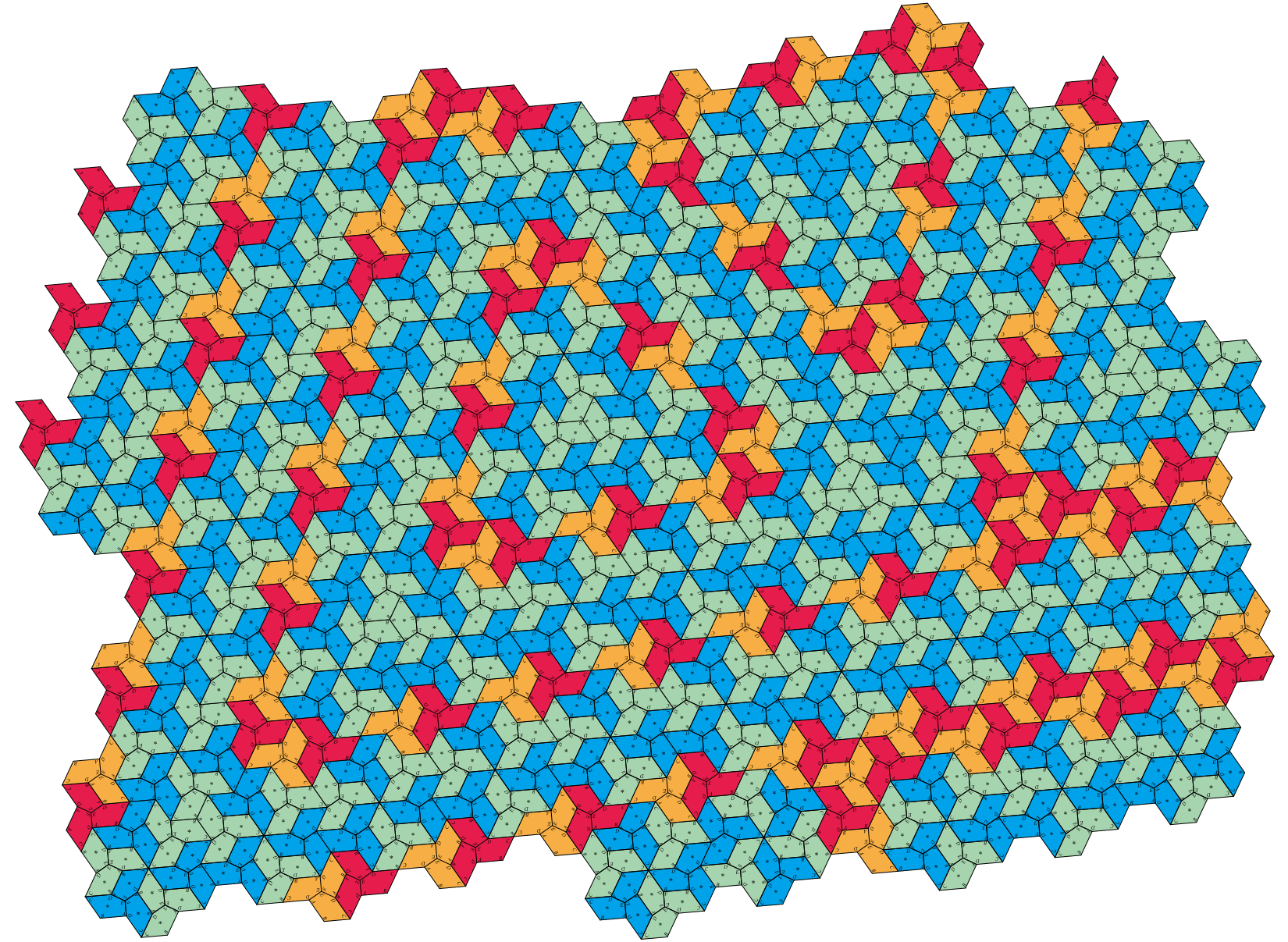} 
  \caption{{\small 
Example of a tiling that is formed by applying the properties of 
Case 2.} 
\label{fig28}
}
\end{figure}

As for Cases 2 and 3 in Figure~\ref{fig22}, it is evident that the unit is formed 
by a hexagonal flowers L1 unit and two reverse windmill units, as shown in 
Figure~\ref{fig24}. By using the units in Figure~\ref{fig24}, the tilings in 
Figure~\ref{fig25} (the properties of Case 3) or Figure~\ref{fig26} (the properties 
of Case 2) can be formed. The tilings as shown in Figures~\ref{fig27} and \ref{fig28} 
can also be formed by applying the property of Case 2. That is, one can arrange 
the width of the reversed windmill units' band or arrange it in a V shape.

\begin{figure}[htbp]
 \centering\includegraphics[width=13.5cm,clip]{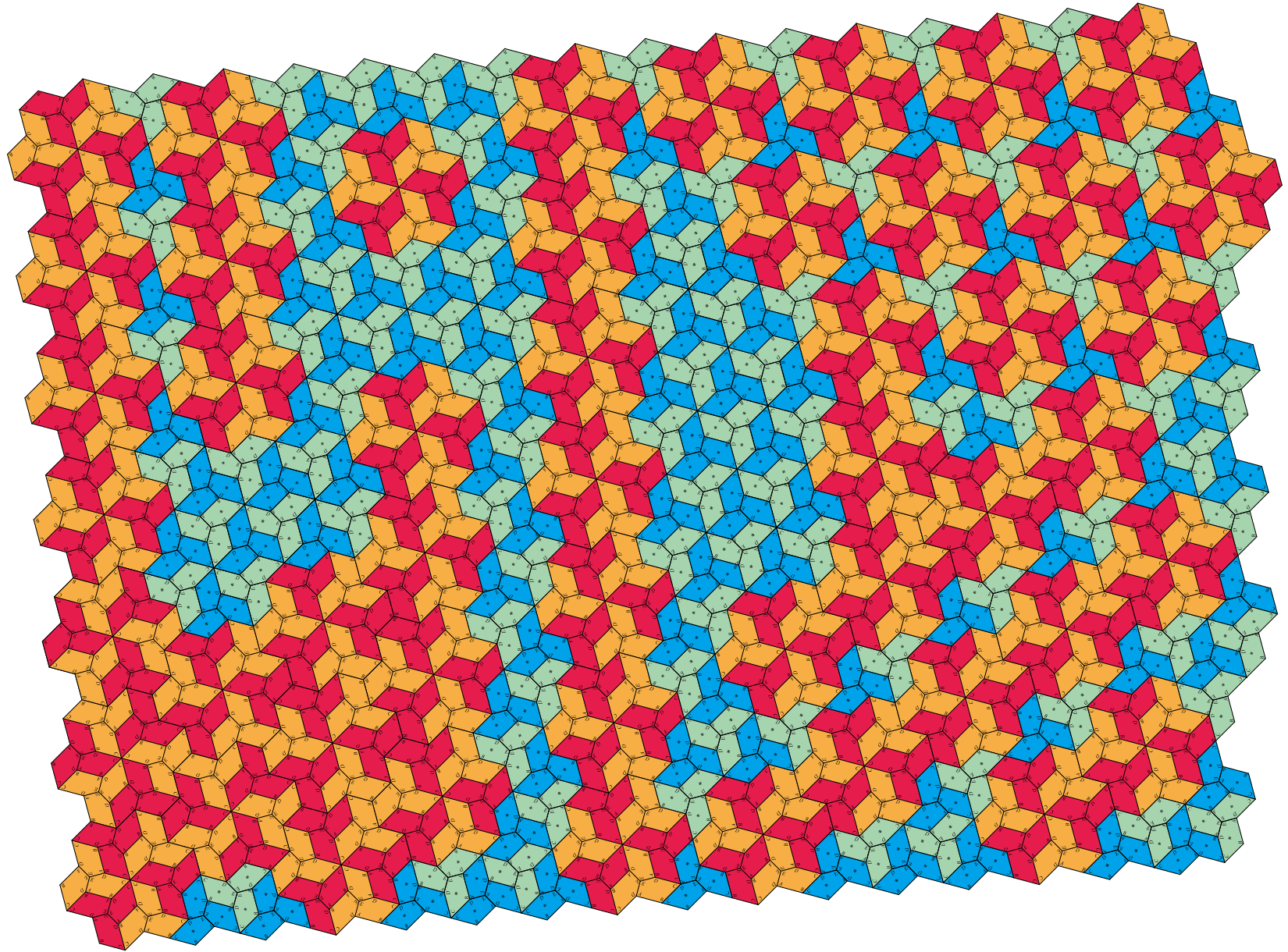} 
  \caption{{\small 
Example of a tiling with the properties of Cases 1, 2 and 3.} 
\label{fig29}
}
\end{figure}

In addition, tilings having the properties of Cases 1, 2, and 3 can also be 
formed. Figure~\ref{fig29} is an example of a tiling with the properties of 
Cases 1, 2, and 3.

Thus, the TH-pentagon admits an infinite variety of periodic tilings and 
nonperiodic tilings by the windmill units.

\begin{figure}[htbp]
 \centering\includegraphics[width=15cm,clip]{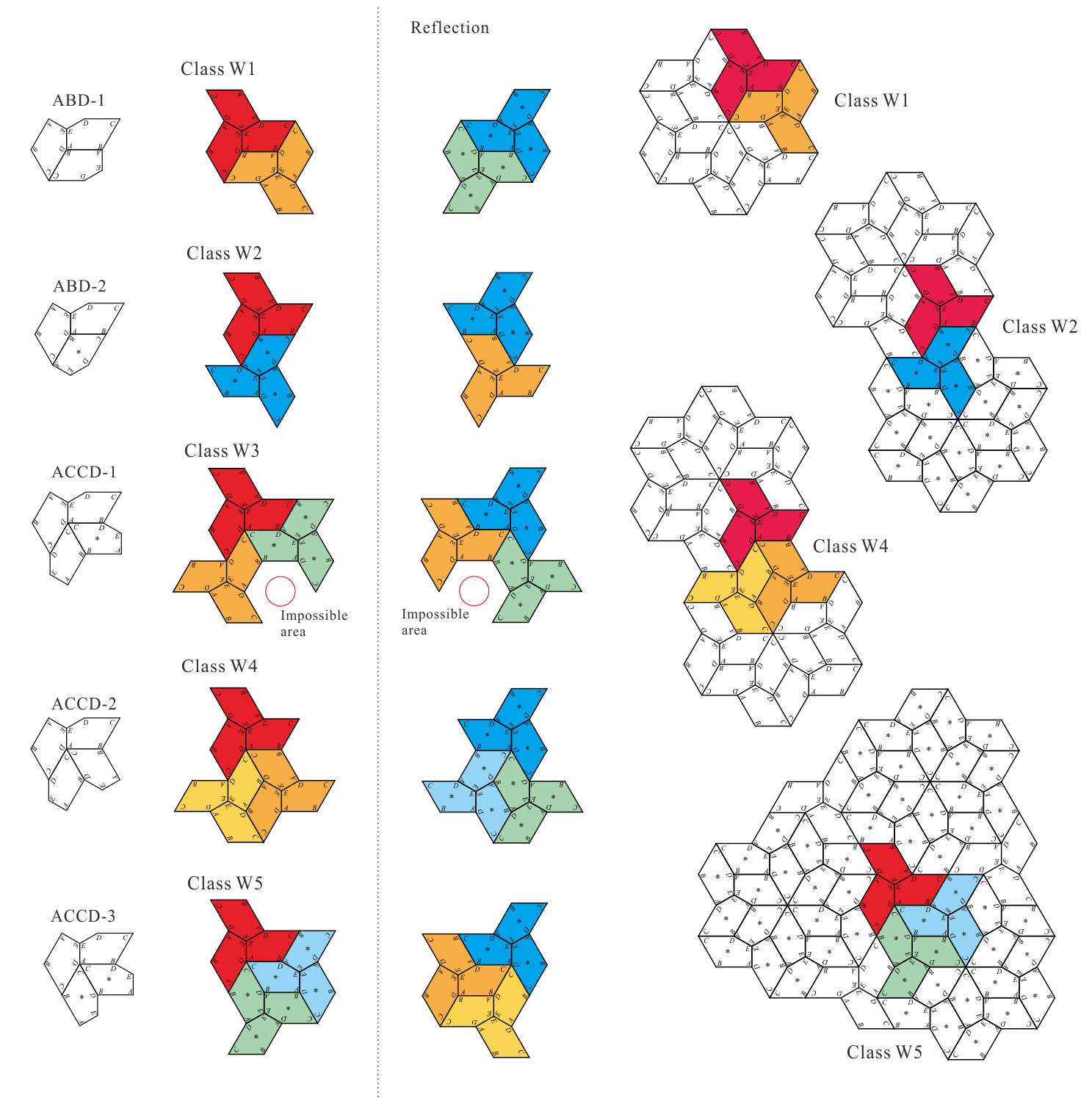} 
  \caption{{\small 
Contiguity method of windmill units.} 
\label{fig30}
}
\end{figure}

Concentration relations of vertices $A$ and $D$ of the TH-pentagon in the 
windmill unit are addressed as follows. When using only the windmill unit, the 
concentration relations of vertices $A$ and $D$ are $A + B + D = 360^ \circ$ 
(ABD-1, ABD-2) or $A + 2C + D = 360^ \circ $ (ACCD-1, ACCD-2, ACCD-3). 
That is, patterns of contiguity method between windmill units are Classes W1, 
W2, W3, W4, and W5 in Figure~\ref{fig30}. Class W1 appears inside the 
HFL1-unit and in a Type 5 tiling. Class W2 appears in a connection of 
AHFL1-unit and PHFL1-unit. Class W3 cannot be used in a tiling. Class W4 a
ppears in a Rice1995-tiling. Class W5 appears in the combination of HFL1-unit 
and the two reverse HFL1-units (see Figure~\ref{fig24}).

\renewcommand{\figurename}{{\small Figure.}}

\renewcommand{\figurename}{{\small Figure.}}

\renewcommand{\figurename}{{\small Figure.}}

\renewcommand{\figurename}{{\small Figure.}}

\renewcommand{\figurename}{{\small Figure.}}

\renewcommand{\figurename}{{\small Figure.}}

\renewcommand{\figurename}{{\small Figure.}}

\renewcommand{\figurename}{{\small Figure.}}

\section{Tilings by only the ship units}
\label{section4}

\begin{figure}[htbp]
 \centering\includegraphics[width=14cm,clip]{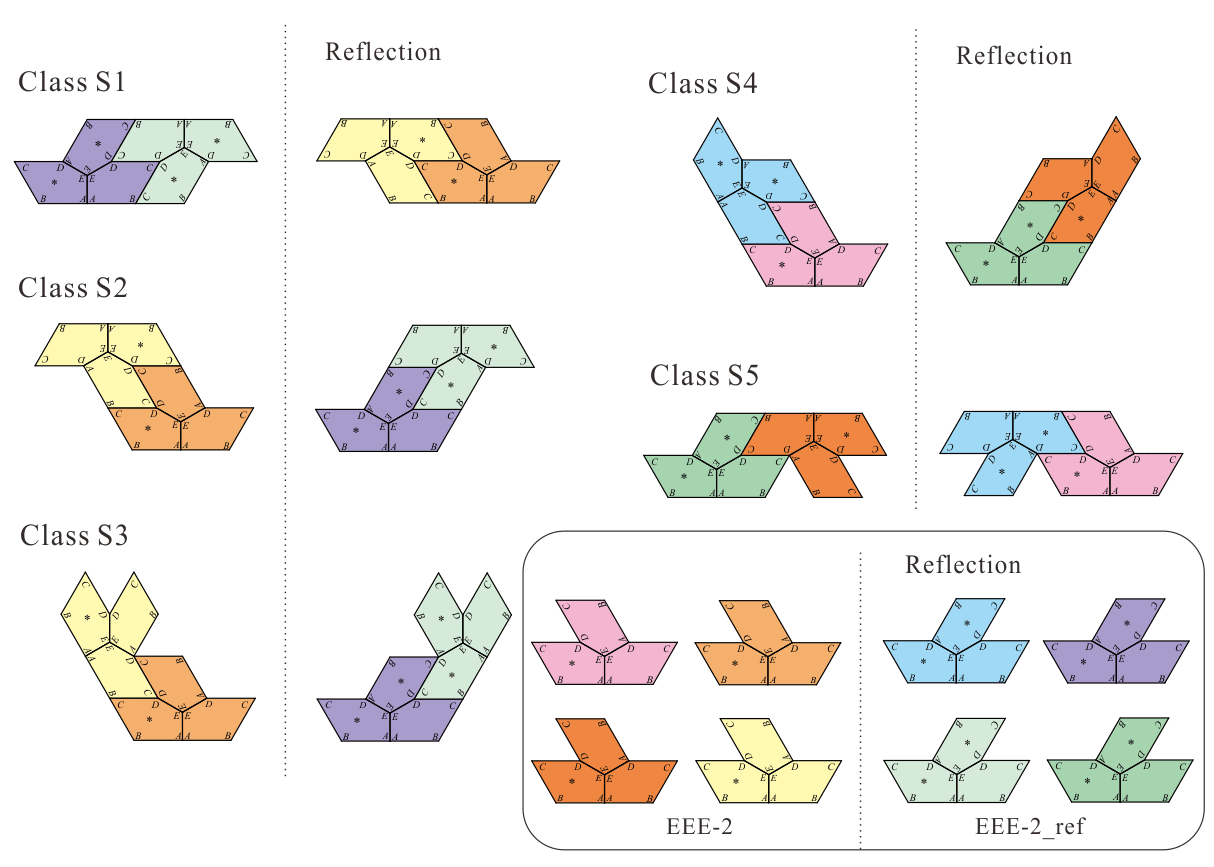} 
  \caption{{\small 
Contiguity method of ship units.} 
\label{fig31}
}
\end{figure}

In this section, tilings by only the ship units are introduced. A concentration 
relation of two vertices $D$ of the TH-pentagon in the ship unit is 
considered next. When using only the ship unit, the concentration relation 
of two vertices $D$ is $C + 2D = 360^ \circ $, and the pattern is CDD 
in Figure~\ref{fig10}. Therefore, the patterns of contiguity method between 
ship units are Classes S1, S2, S3, S4, and S5 in Figure~\ref{fig31}.

\begin{figure}[htbp]
 \centering\includegraphics[width=15cm,clip]{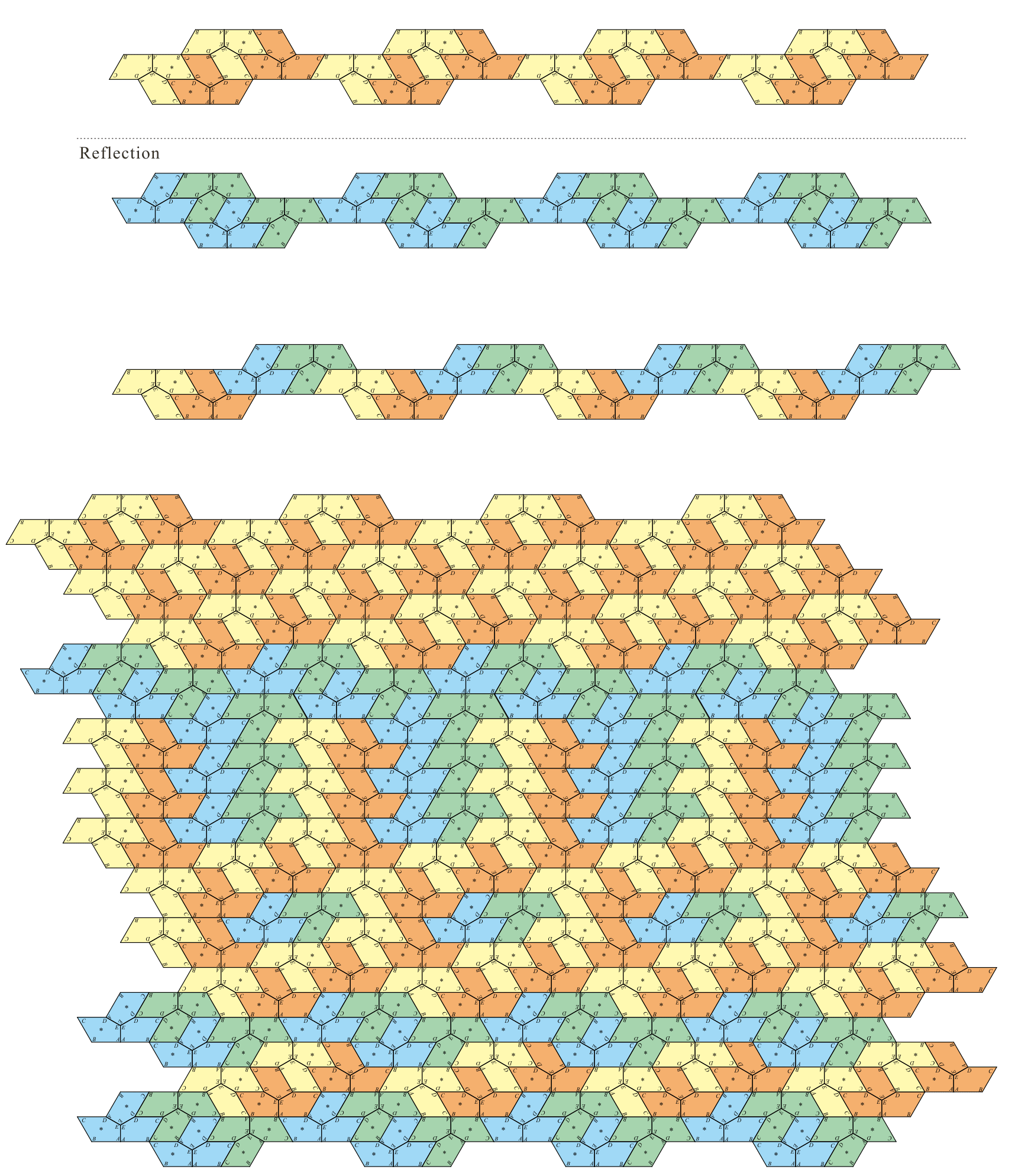} 
  \caption{{\small 
Tiling with Class S1.} 
\label{fig32}
}
\end{figure}

\begin{figure}[htbp]
 \centering\includegraphics[width=13.5cm,clip]{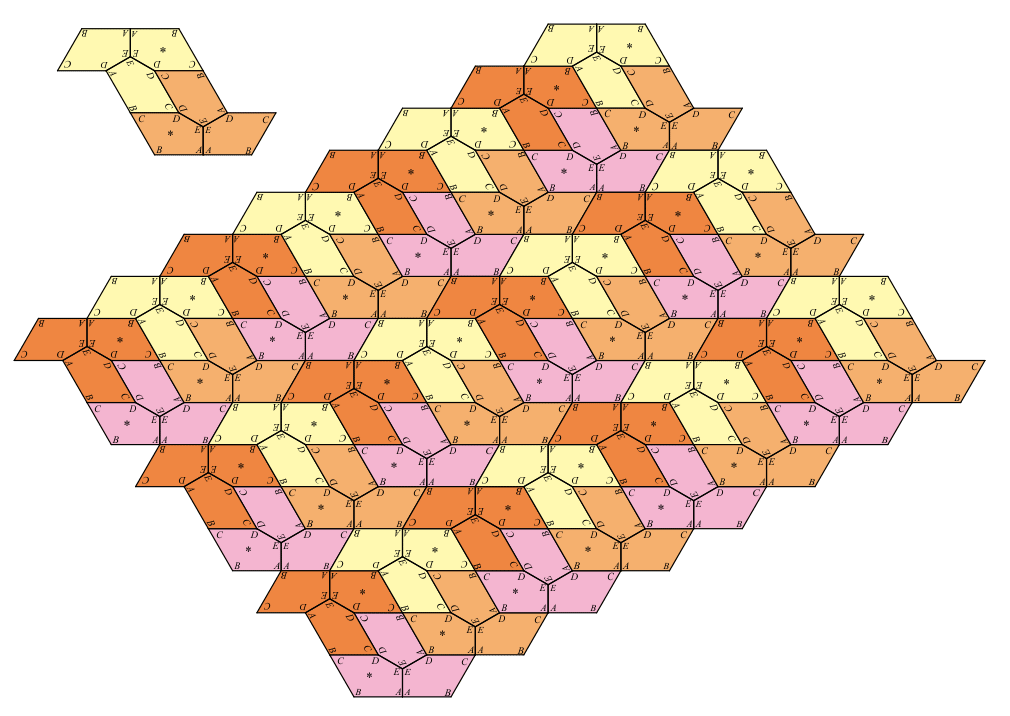} 
  \caption{{\small 
Tiling with Class S2.} 
\label{fig33}
}
\end{figure}

\begin{figure}[htbp]
 \centering\includegraphics[width=15.5cm,clip]{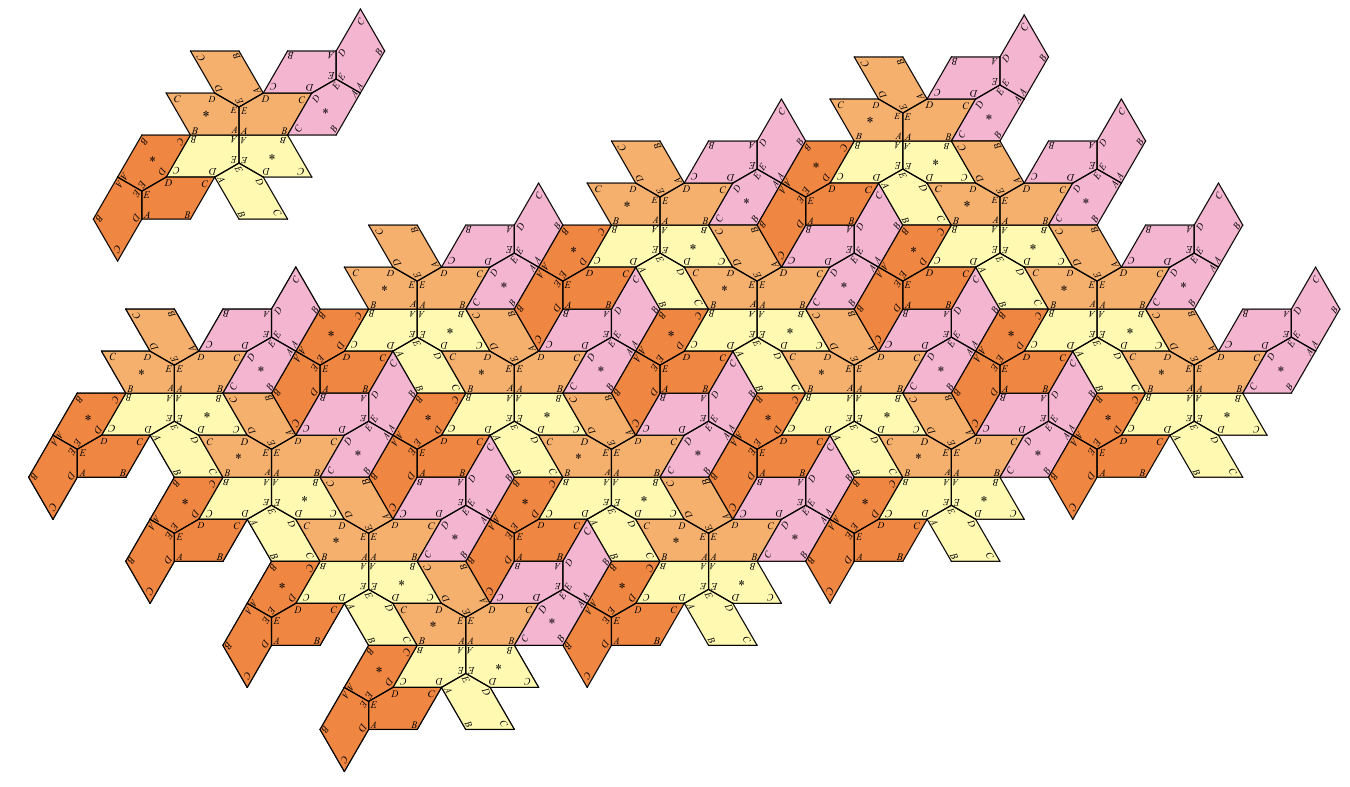} 
  \caption{{\small 
Tiling with Class S3.} 
\label{fig34}
}
\end{figure}

\begin{figure}[htbp]
 \centering\includegraphics[width=15cm,clip]{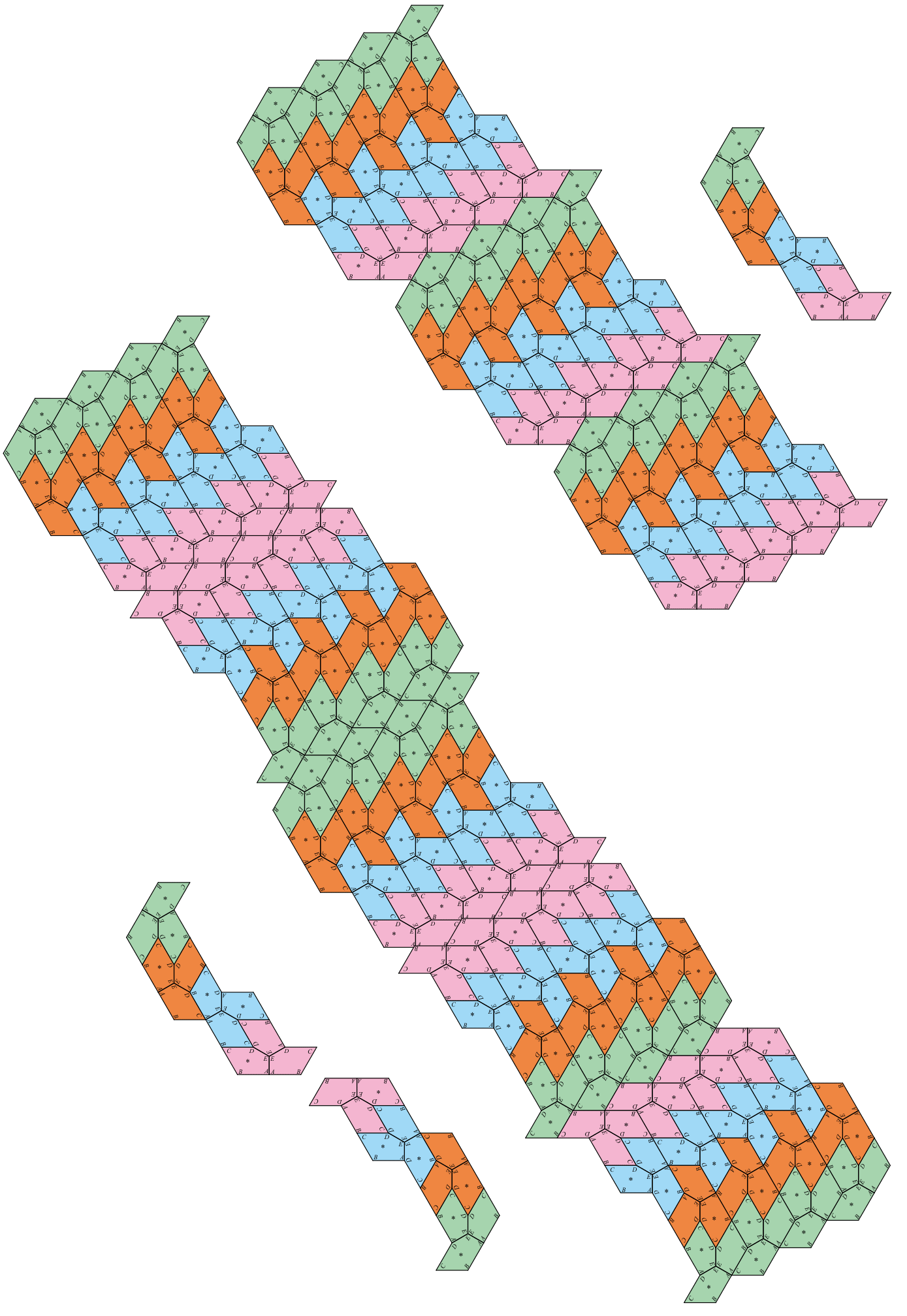} 
  \caption{{\small 
Tiling with Class S4.} 
\label{fig35}
}
\end{figure}

\begin{figure}[htbp]
 \centering\includegraphics[width=10.5cm,clip]{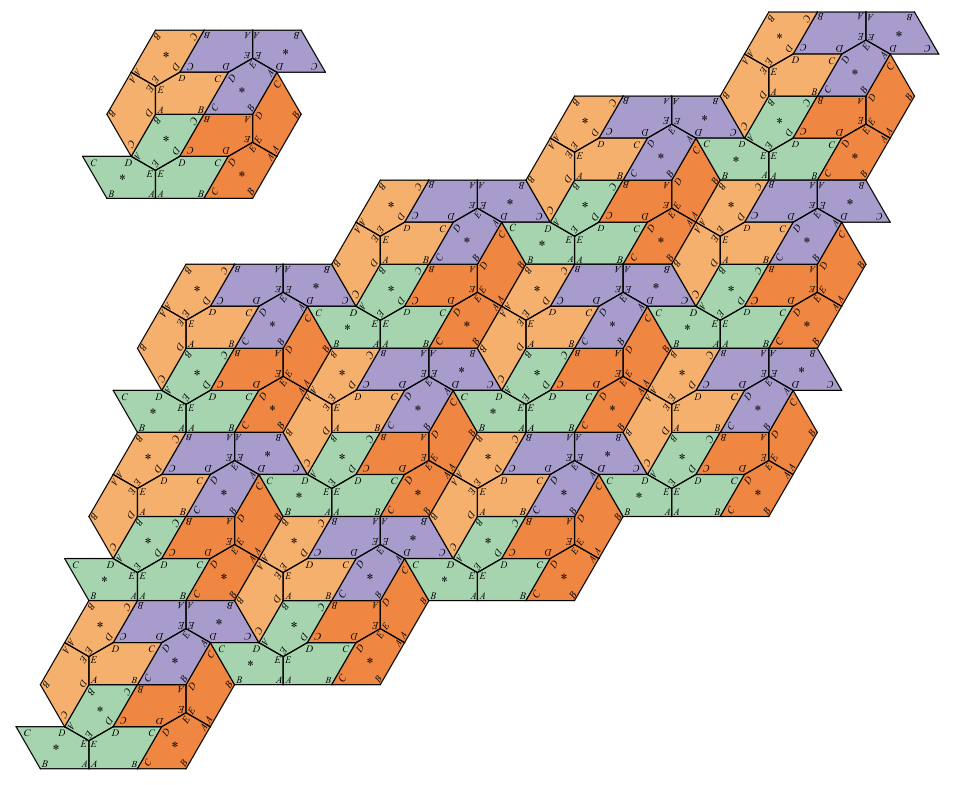} 
  \caption{{\small 
Tiling with Class S4.} 
\label{fig36}
}
\end{figure}

\begin{figure}[htbp]
 \centering\includegraphics[width=13.5cm,clip]{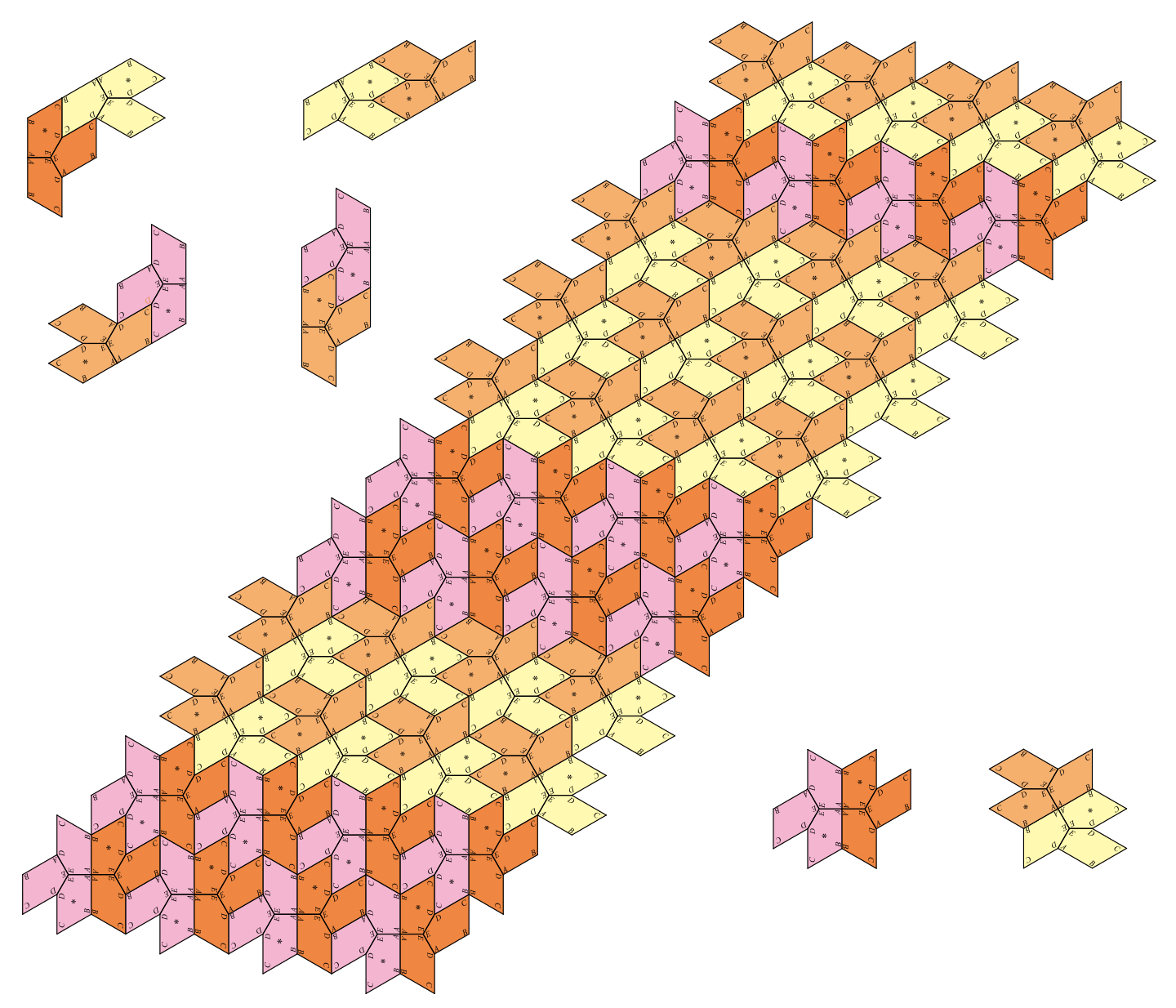} 
  \caption{{\small 
Tiling with Classes S1 and S3.} 
\label{fig37}
}
\end{figure}

\begin{figure}[htbp]
 \centering\includegraphics[width=14cm,clip]{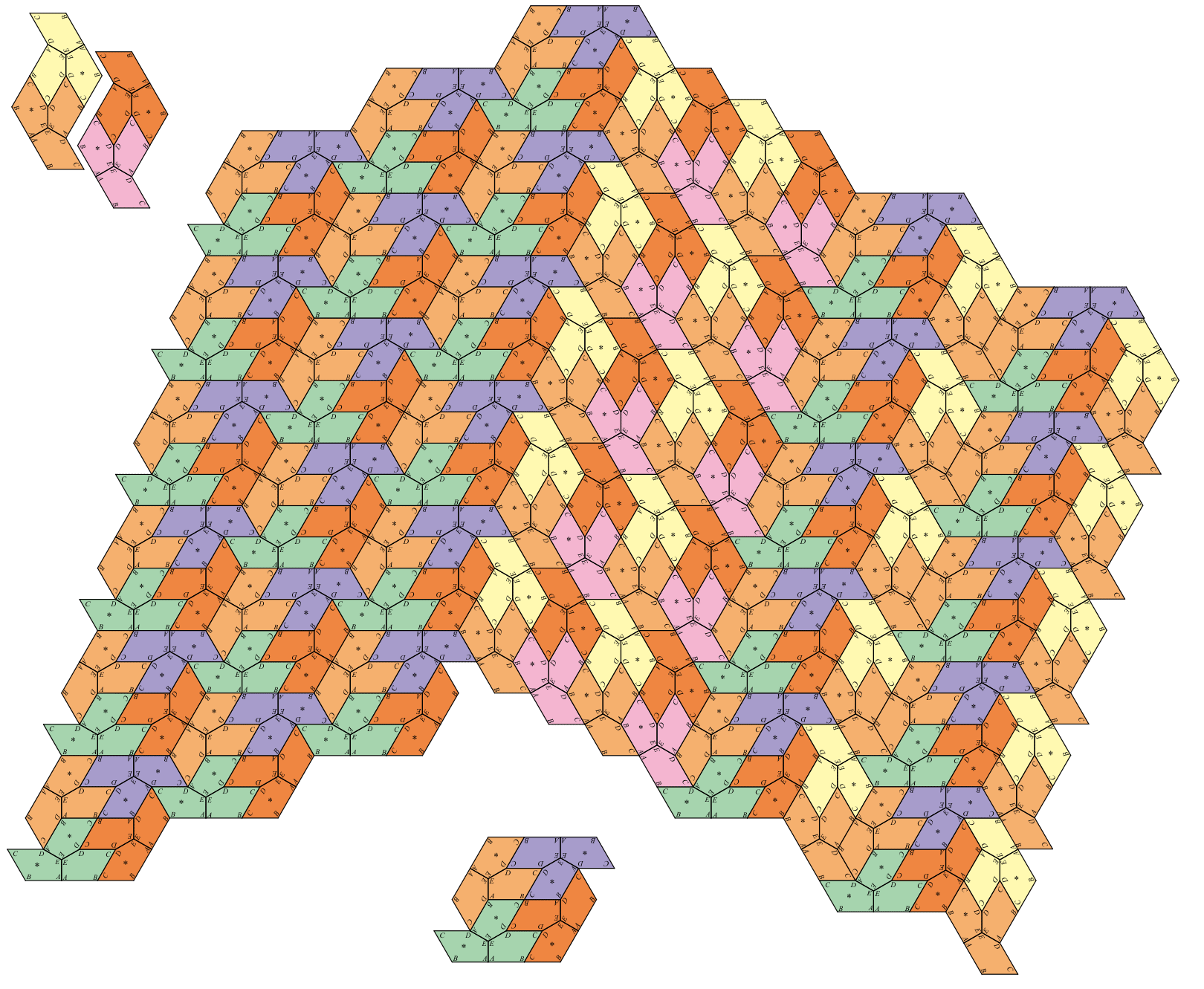} 
  \caption{{\small 
Tiling with Classes S2 and S4.} 
\label{fig38}
}
\end{figure}

The tiling in Figure~\ref{fig32} is a tiling with Class S1. The tiling of Figure~\ref{fig32} 
has the property of freely arranging three kinds of objects like a band in 
one-dimension. The tiling in Figure~\ref{fig33} is a tiling with Class S2, and the 
tiling in Figure~\ref{fig34} is a tiling with Class S3. The tilings in Figures~\ref{fig35} and 
\ref{fig36} are tilings with Class S4. Then, the tiling of Figure~\ref{fig35} has the property 
of freely arranging base parts in one-dimension. There is no tiling with 
only Class S5. In addition, the tiling in Figure~\ref{fig37} is a tiling with Classes 
S1 and S3. The tiling in Figure~\ref{fig38} is a tiling with Classes S2 and S4. The 
tilings of Figures~\ref{fig37} and \ref{fig38} have the property that base parts can 
combine freely in one dimension.

 Those tilings are the tilings of only the ship unit the authors have 
currently found. (In fact, although there are other tilings, they are shown 
in the next section by their properties.)

 Thus, the TH-pentagon admits an infinite variety of periodic tilings and 
nonperiodic tilings by the ship units.

\renewcommand{\figurename}{{\small Figure.}}

\renewcommand{\figurename}{{\small Figure.}}

\renewcommand{\figurename}{{\small Figure.}}

\renewcommand{\figurename}{{\small Figure.}}

\renewcommand{\figurename}{{\small Figure.}}

\renewcommand{\figurename}{{\small Figure.}}

\renewcommand{\figurename}{{\small Figure.}}

\renewcommand{\figurename}{{\small Figure.}}

\renewcommand{\figurename}{{\small Figure.}}
\begin{figure}[htbp]
 \centering\includegraphics[width=14.5cm,clip]{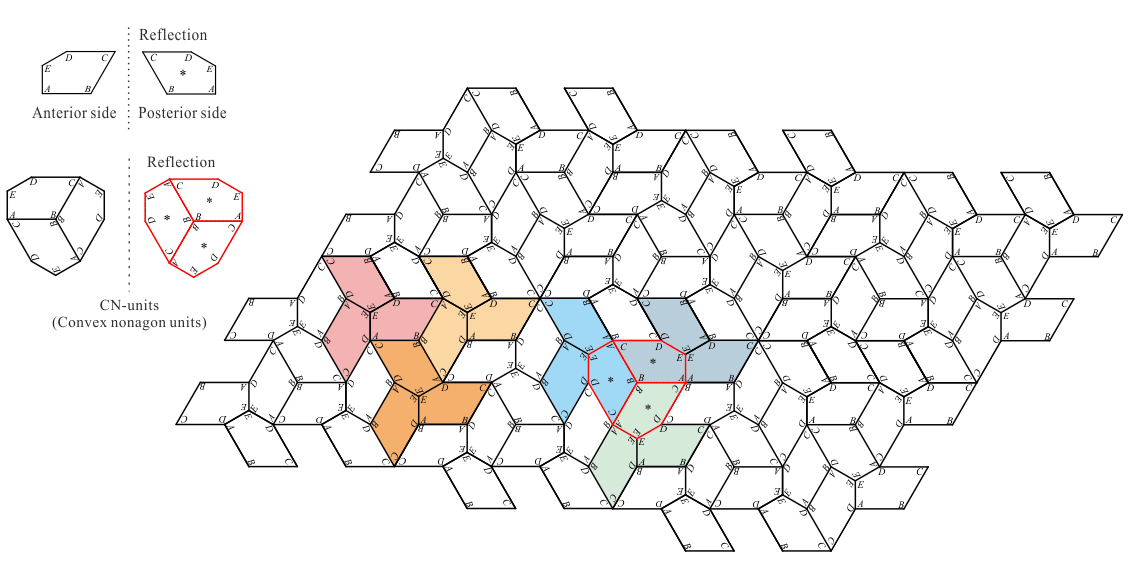} 
  \caption{{\small 
Reversing of a CN-unit in Rice1995-tiling.} 
\label{fig39}
}
\end{figure}

\section{Tilings by the windmill units and the ship units}
\label{section5}

In this section, the authors introduce the result that the TH-pentagon 
admits infinite variety of periodic tilings and nonperiodic tilings by the 
windmill units and the ship units.

\subsection{Rice1995-tiling with ACN-units and PCN-units }
\label{subsection5_1}

The Rice1995-tiling in Figure~\ref{fig4} is a Rice1995-tiling by the anterior 
TH-pentagon, which is considered as tiling by only the anterior windmill 
unit. The Rice1995-tiling in Figure~\ref{fig39} contains a PCN-unit. Apparently, it 
is equal to replacing a windmill unit with a ship unit (or vice versa) when 
a CN-unit is reversed. Therefore, the Rice1995-tilings with ACN-units and 
PCN-units as shown in Figures~\ref{fig7} and \ref{fig39} are tilings by the windmill units and 
the ship units.

\subsection{Consideration of similar hexagonal flower units }
\label{subsection5_2}

\begin{figure}[htbp]
 \centering\includegraphics[width=14cm,clip]{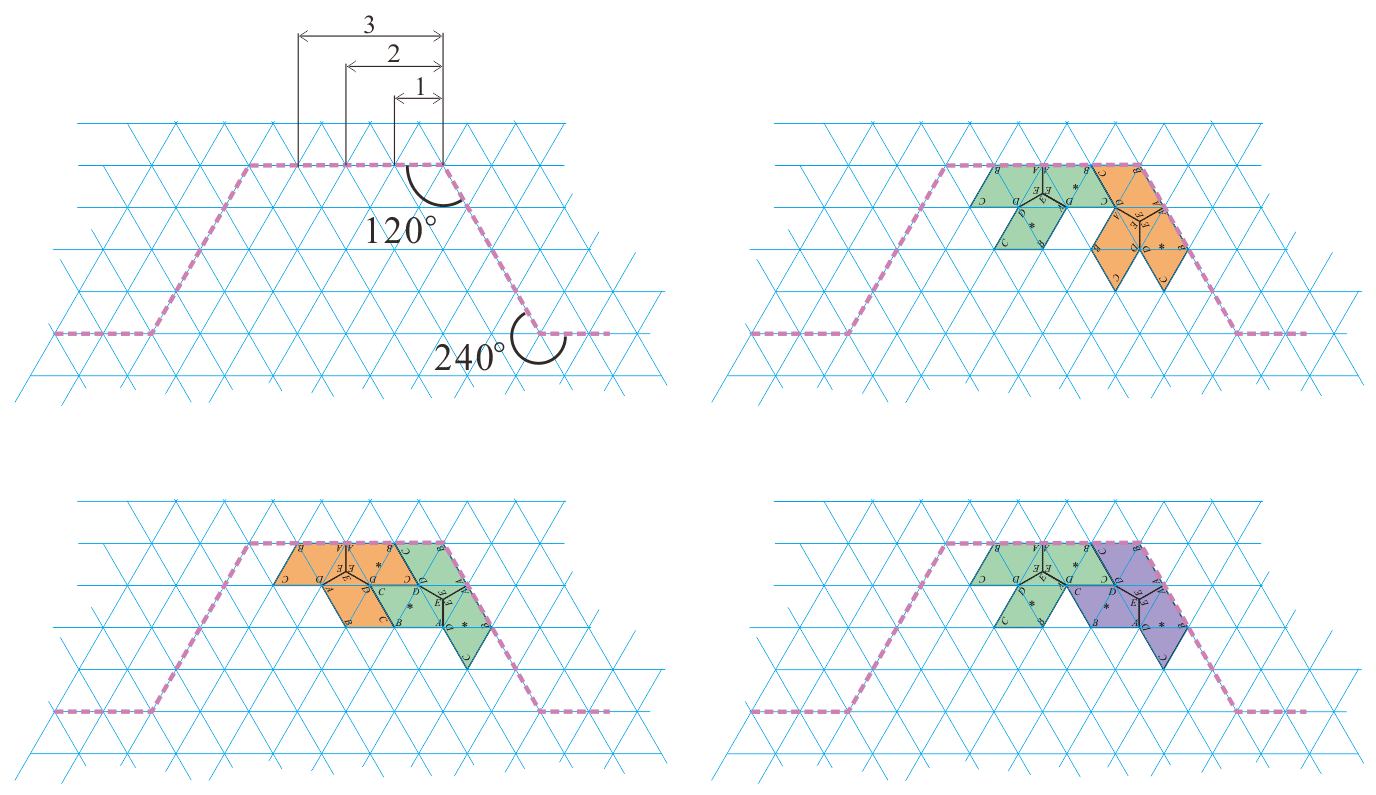} 
  \caption{{\small 
Three patterns that use the two ship units.} 
\label{fig40}
}
\end{figure}

In Section~\ref{section3}, the hexagonal flower unit (i.e., SHL1-unit) formed by six 
windmill units was presented. As shown in the next subsections, by using the 
windmill units and ship units, it is possible to create a similar hexagonal 
flower unit with the side length of two or three. On the other hand, the 
authors have confirmed that there are no hexagonal flower units with the 
side length of four or more. The reason is briefly described below. Focusing 
on the convex portion of the hexagonal flower unit, one obtains combinations 
of windmill units and ship units capable of forming the convex portion. In 
addition, the vertex of $120^ \circ $ in the convex portion will use $B = 
120^ \circ $ or $C + C = 120^ \circ $. As a result, the combinations of 
which the side length is three or more and will generate the tiling are 
three patterns (not distinguishing reflections and rotations) that use the 
two ship units in Figure~\ref{fig40}. By using these three patterns, the authors 
confirmed that there are no hexagonal flower units with the side length of 
five or more. In addition, when the side length is four, it is possible to 
form a part, but when trying to form hexagonal flower units with side length 
of four, an impossible area appears\footnote{ By also a program that used 
Algorithm X, the authors confirmed that there is no hexagonal flower unit with 
side length four~\cite{Hiroi_a, Hiroi_b, Knuth_2000, wiki-10, wiki-11, wiki-12}.}.
Therefore, it is impossible to form hexagonal flower units with side 
lengths of four or more by using TH-pentagons.

\renewcommand{\figurename}{{\small Figure.}}

\subsection{Hexagonal flower L2 unit and tilings}
\label{subsection5_3}

\begin{figure}[htbp]
 \centering\includegraphics[width=15cm,clip]{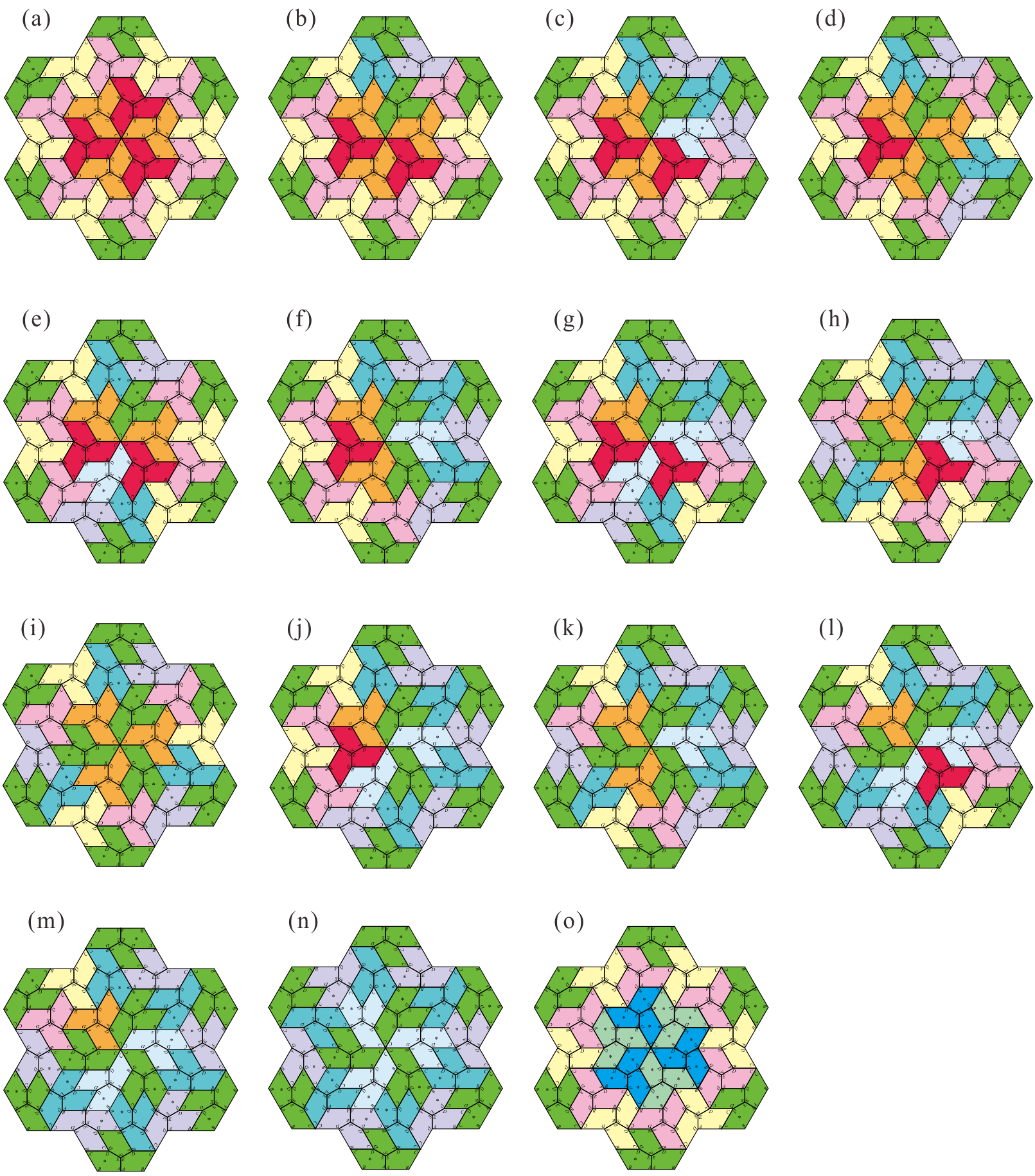} 
  \caption{{\small 
HFL2-units of 15 unique patterns.} 
\label{fig41}
}
\end{figure}

Hereafter, the shape that is formed by 72 TH-pentagons in Figure~\ref{fig41} is 
referred to as \textit{a hexagonal flowers L2 unit} (HFL2-unit). The length of 
the sides of HFL2-unit is twice the length of the sides of HFL1-unit. 
The arrangements of the internal TH-pentagons (or windmill units and 
ship units) that form HFL2-unit are 15 patterns in Figure~\ref{fig41}, 
not distinguishing reflections and rotations (i.e., the number of unique 
patterns of HFL2-unit is 15). Note that it was confirmed by the 
program that there are 15 unique patterns.

\begin{figure}[htbp]
 \centering\includegraphics[width=12cm,clip]{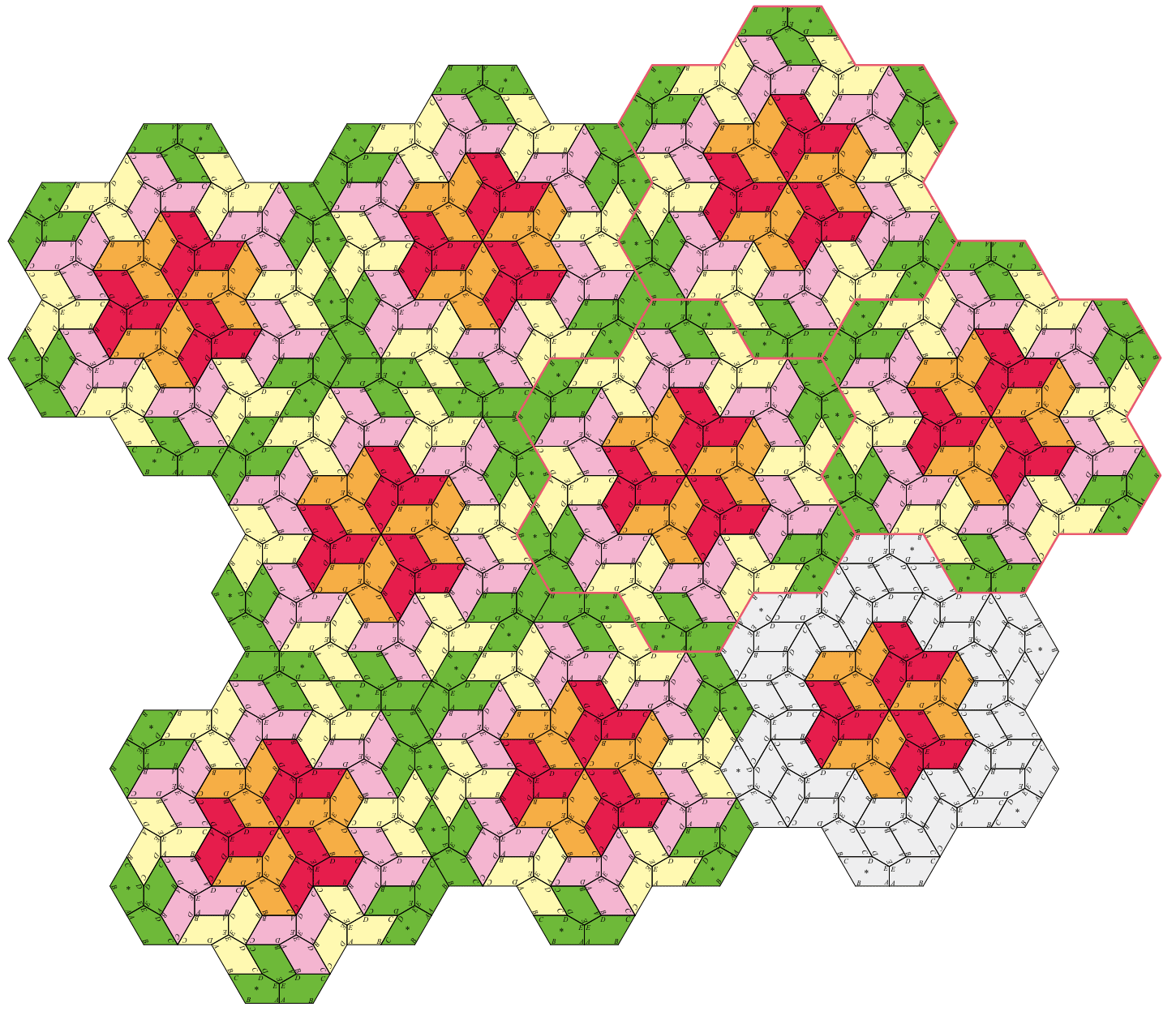} 
  \caption{{\small 
Example of tiling by HFL2-units.} 
\label{fig42}
}
\end{figure}

\begin{figure}[htbp]
 \centering\includegraphics[width=12cm,clip]{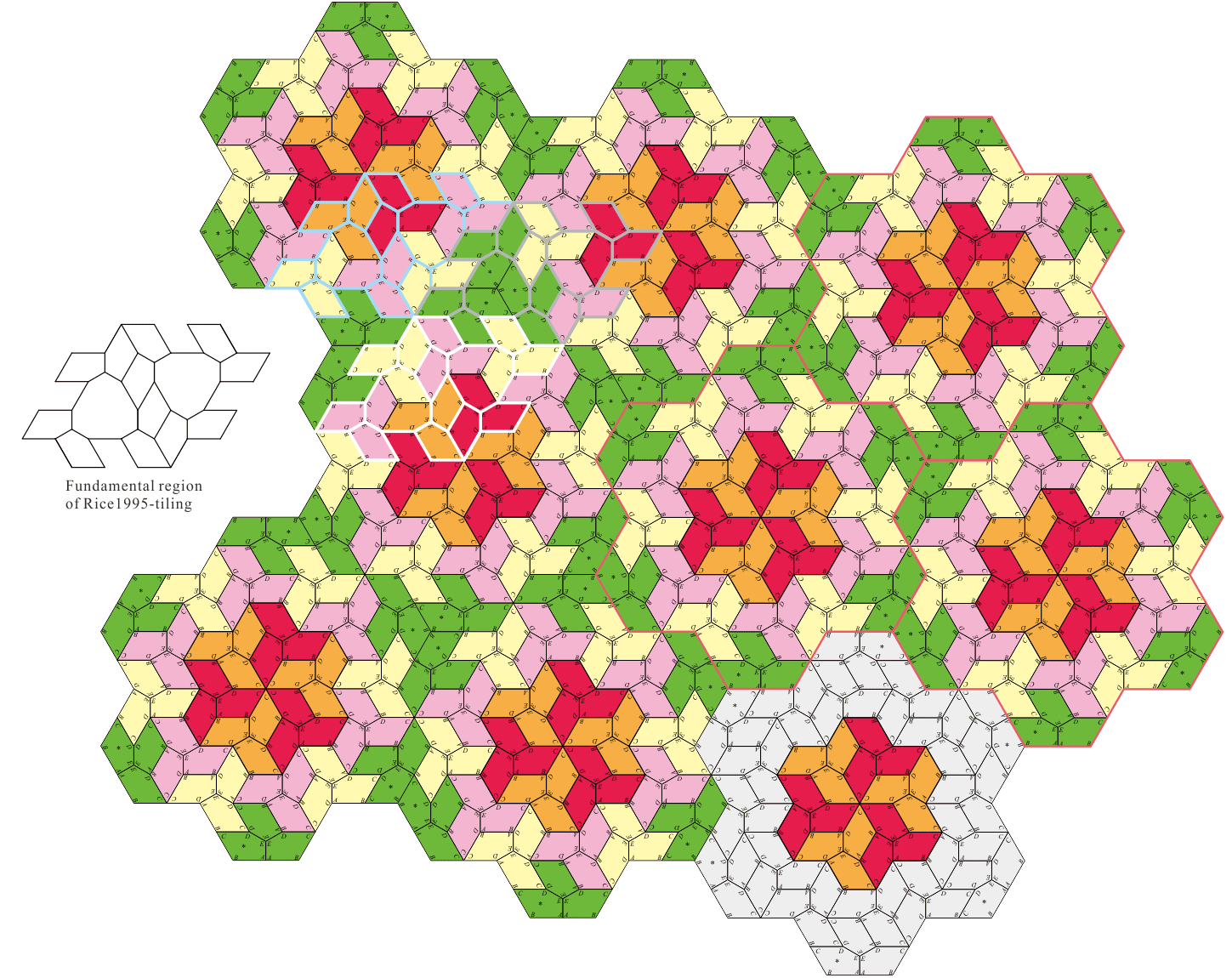} 
  \caption{{\small 
Example of tiling by HFL2-units.} 
\label{fig43}
}
\end{figure}

\begin{figure}[htbp]
 \centering\includegraphics[width=15cm,clip]{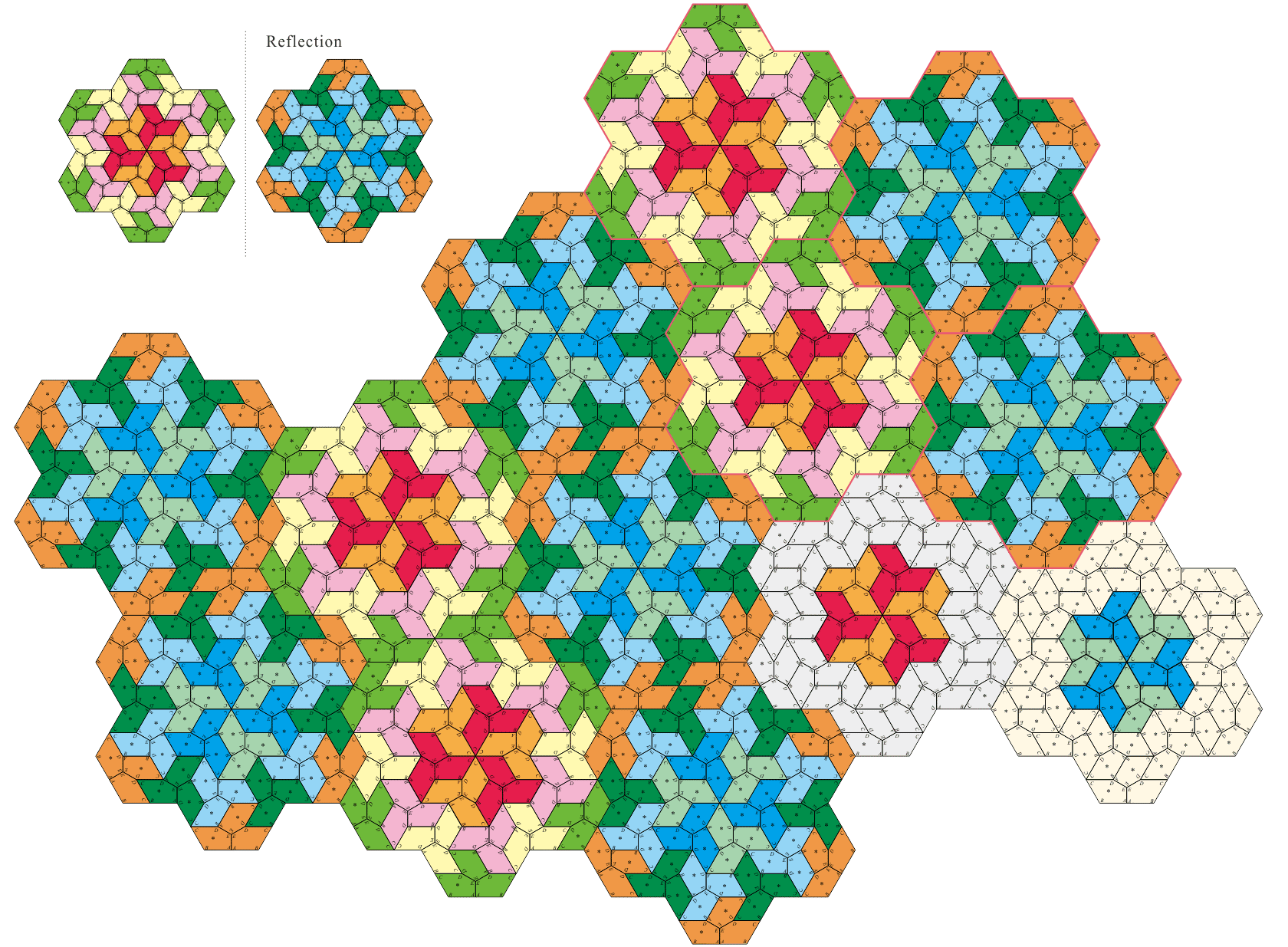} 
  \caption{{\small 
Example of tiling by HFL2-units.} 
\label{fig44}
}
\end{figure}

\begin{figure}[htbp]
 \centering\includegraphics[width=15cm,clip]{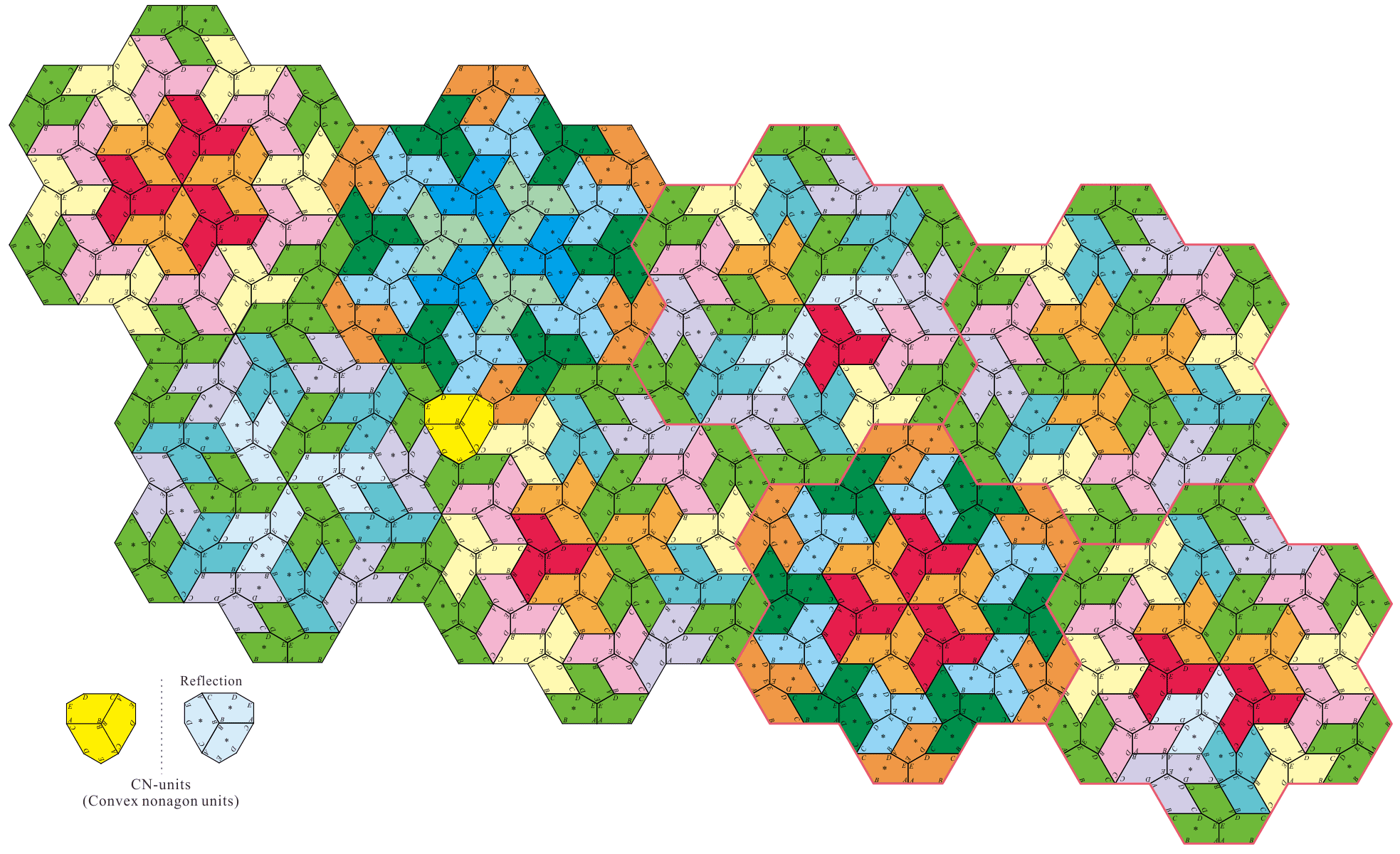} 
  \caption{{\small 
Example of tiling by HFL2-units.} 
\label{fig45}
}
\end{figure}

Like HFL1-units, the HFL2-units can generate tilings by combining 
HFL2-units. The contiguity methods are three patterns as shown in 
Figures~\ref{fig42}, \ref{fig43}, and \ref{fig44}. Note that the tilings in 
Figures~\ref{fig42}, \ref{fig43}, and \ref{fig44} are formed by an 
HFL2-unit of Figure~\ref{fig41}(a) and its reflected image. (As for tilings in 
Figure~\ref{fig56} in the Appendix, they are cases of an HFL2-unit of 
Figure~\ref{fig41}(n).) Of course, the parts of CN-units in tilings can 
be reversed freely, and different patterns of HFL2-units which 
contain reflections and rotations can be freely connectable in the 
one tiling (see Figure~\ref{fig45}).

HFL2-units in Figure~\ref{fig41}(a) and (o) contain a concentric HFL1-unit, and 
contain 18 windmill units and six ship units. HFL2-units in Figure~\ref{fig41}(a)--(n) 
contain six CN-units. HFL2-units in Figure~\ref{fig41}(n) are formed by 24 ship 
units. Therefore, a tiling by only the HFL2-units in Figure~\ref{fig41}(n) (see 
Figure~\ref{fig56} in Appendix) is a tiling with CN-units by only the ship 
units\footnote{ If all the CN-units in the Rice1995-tiling in Figure~\ref{fig4} is 
reflected, it is a Rice1995-tiling by only the ship units. This tiling is 
equal to the tiling in Figure~\ref{fig56}(b) in Appendix. That is, a tiling of 
HFL2-units according to the contiguity method in Figure~\ref{fig43} and 
Figure~\ref{fig56}(b) is a Rice1995-tiling. }. 

The outline of the fundamental region of the Rice1995-tiling appears as a 
unit in the tiling in Figure~\ref{fig43}. Therefore, as shown in Figure~\ref{fig46} or 
Figures~\ref{fig57} and ~\ref{fig58} in the Appendix, tilings with HFL1-units and 
tilings with HFL2-units are connectable.

\renewcommand{\figurename}{{\small Figure.}}

\renewcommand{\figurename}{{\small Figure.}}

\renewcommand{\figurename}{{\small Figure.}}

\renewcommand{\figurename}{{\small Figure.}}

\renewcommand{\figurename}{{\small Figure.}}

\renewcommand{\figurename}{{\small Figure.}}
\begin{figure}[htbp]
 \centering\includegraphics[width=15cm,clip]{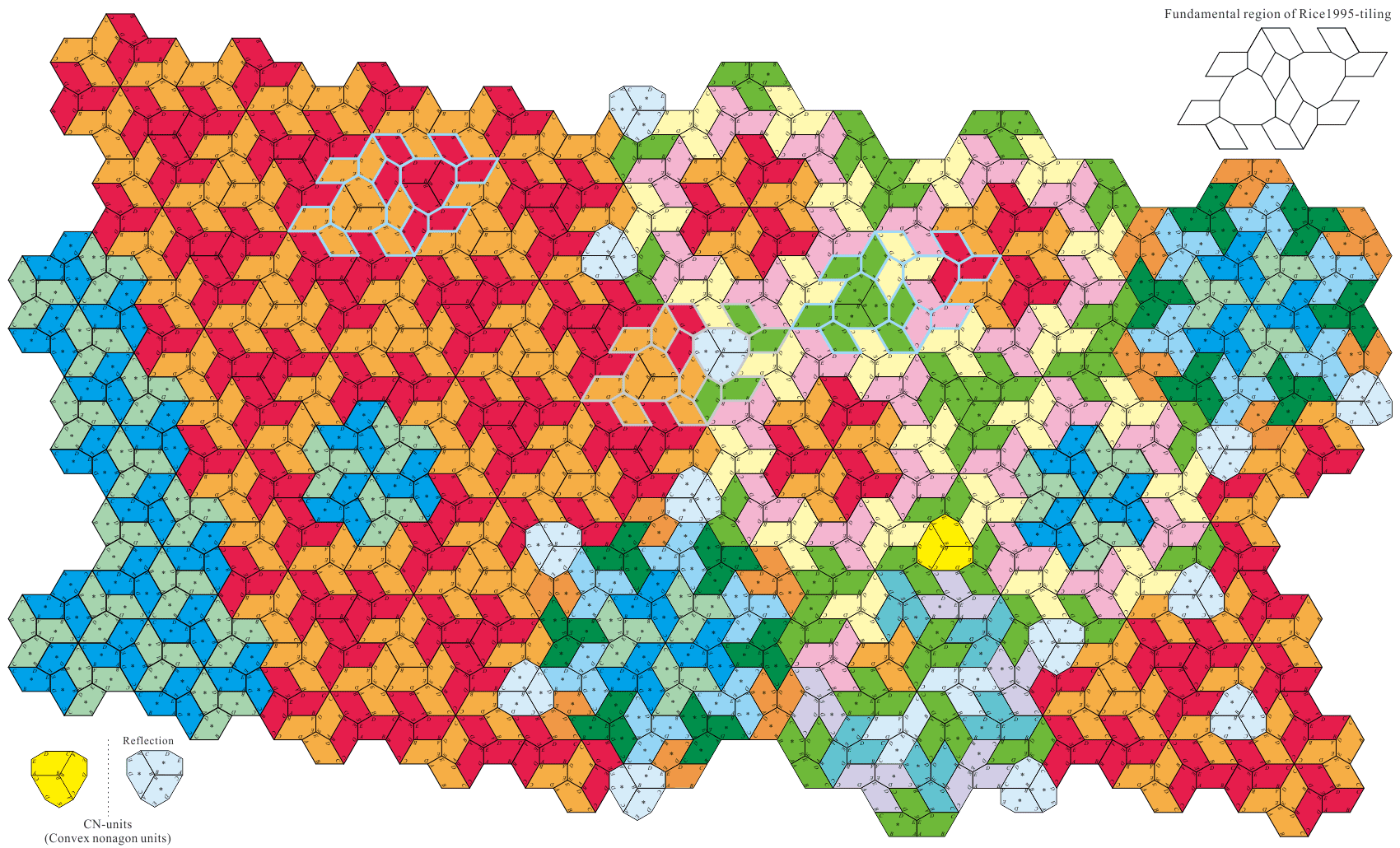} 
  \caption{{\small 
Example of tiling by HFL1-units and HFL2-units.} 
\label{fig46}
}
\end{figure}

\subsection{Hexagonal flower L3 unit and tilings}
\label{subsection5_4}

\begin{figure}[htbp]
 \centering\includegraphics[width=14.5cm,clip]{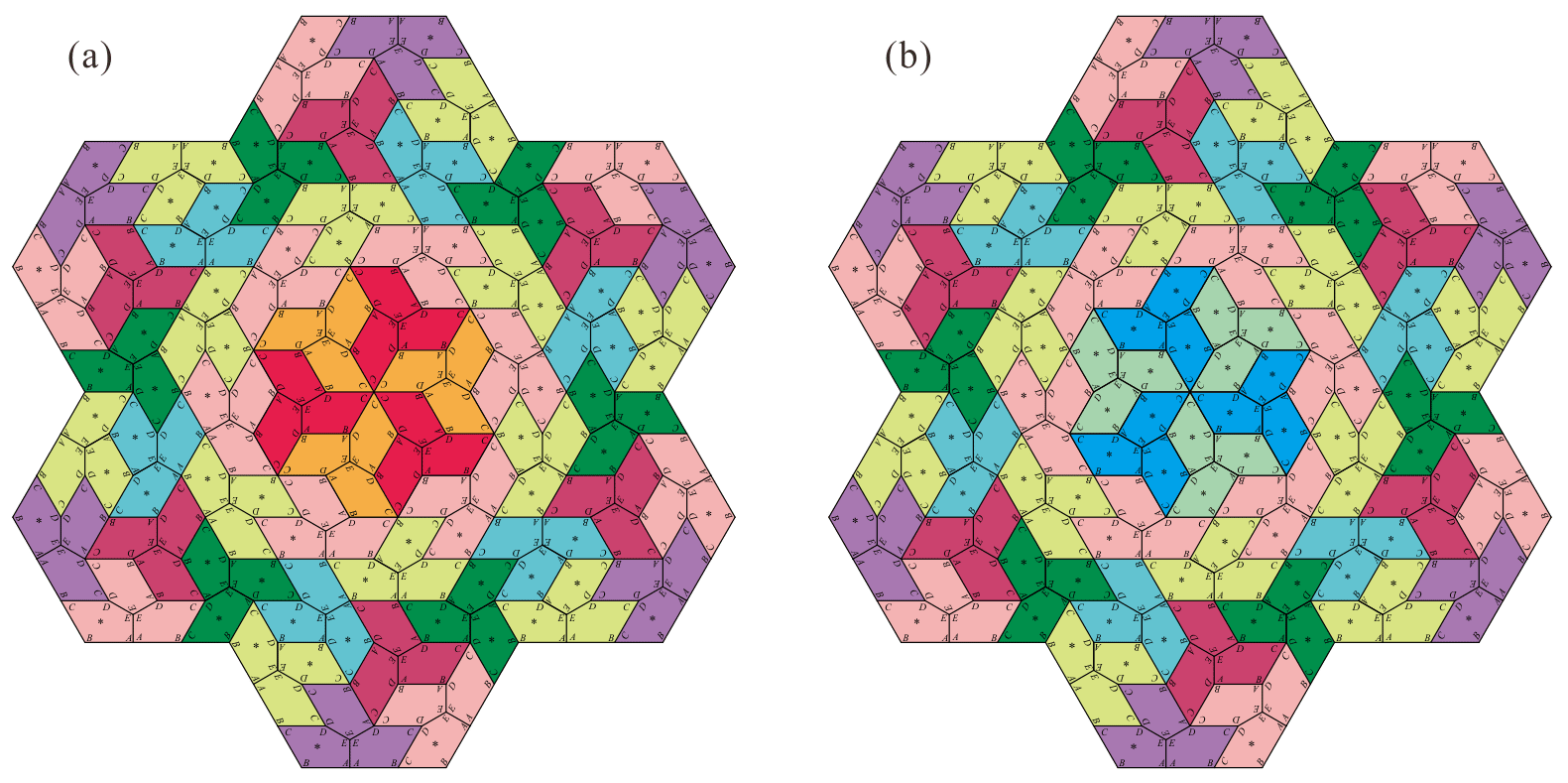} 
  \caption{{\small 
HFL3-units of two unique patterns.} 
\label{fig47}
}
\end{figure}

Hereafter, the shape that is formed by 162 TH-pentagons in Figure~\ref{fig47} is 
referred to as \textit{a hexagonal flowers L3 unit} (HFL3-unit)\footnote{ The 
arrangement of the internal TH-pentagons of HFL3-units was discovered by 
Toshihiro Shirakawa.}. The length of the sides of an HFL3-unit is three times 
the length of sides of an HFL1-unit. The arrangements of the internal 
TH-pentagons (or windmill units and ship units) that form HFL3-unit are two 
patterns in Figure~\ref{fig47}, not distinguishing reflections and rotations 
(i.e., the number of unique patterns of HFL3-unit is two). Note that it was 
confirmed by the program that there are two unique patterns. HFL3-units 
in Figure~\ref{fig47} contain a concentric HFL1-unit, and contain 18 
windmill units and 36 ship units.

\begin{figure}[htbp]
 \centering\includegraphics[width=13.5cm,clip]{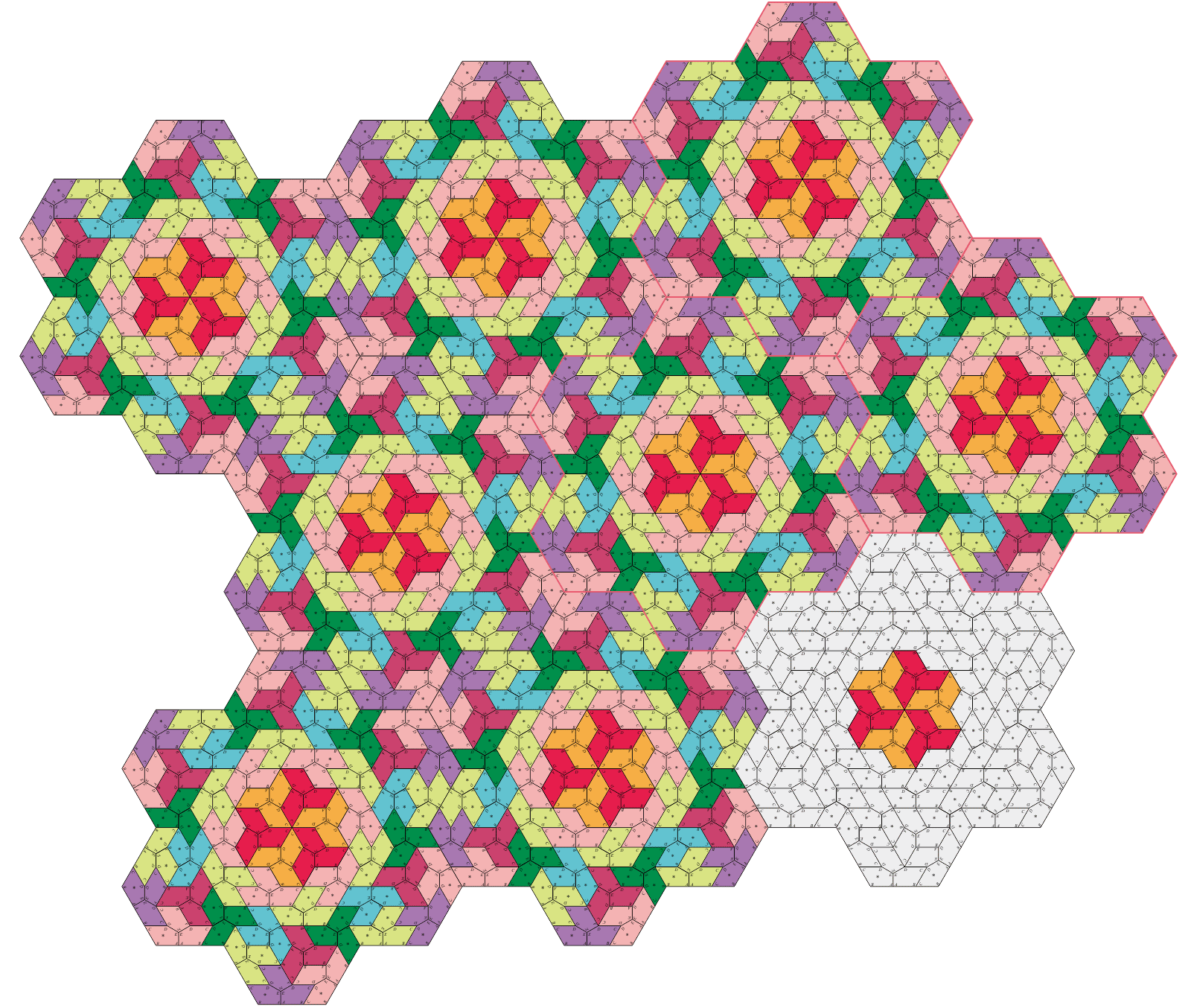} 
  \caption{{\small 
Example of tiling by HFL3-units.} 
\label{fig48}
}
\end{figure}

\begin{figure}[htbp]
 \centering\includegraphics[width=12cm,clip]{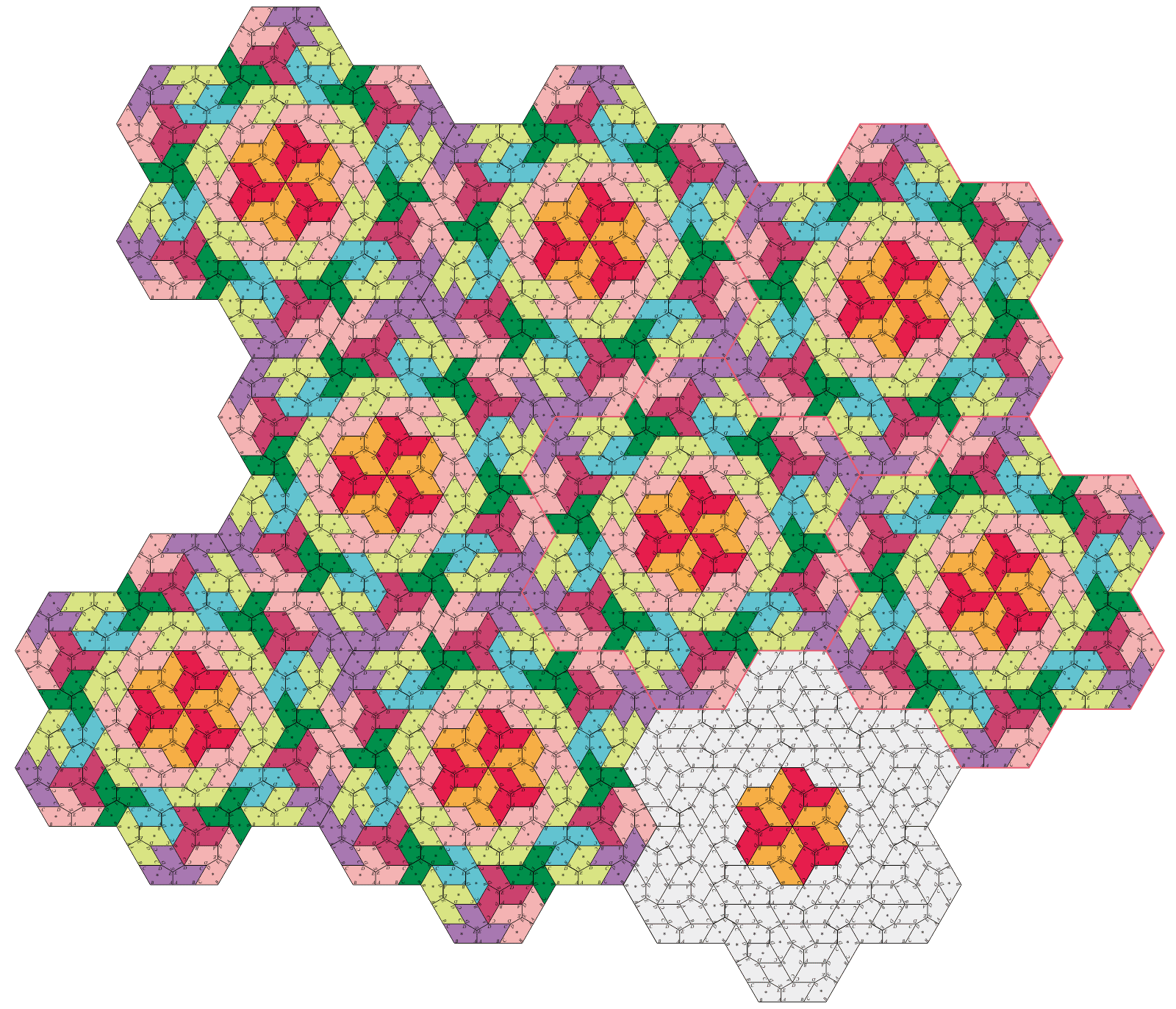} 
  \caption{{\small 
Example of tiling by HFL3-units.} 
\label{fig49}
}
\end{figure}

\begin{figure}[htbp]
 \centering\includegraphics[width=13.5cm,clip]{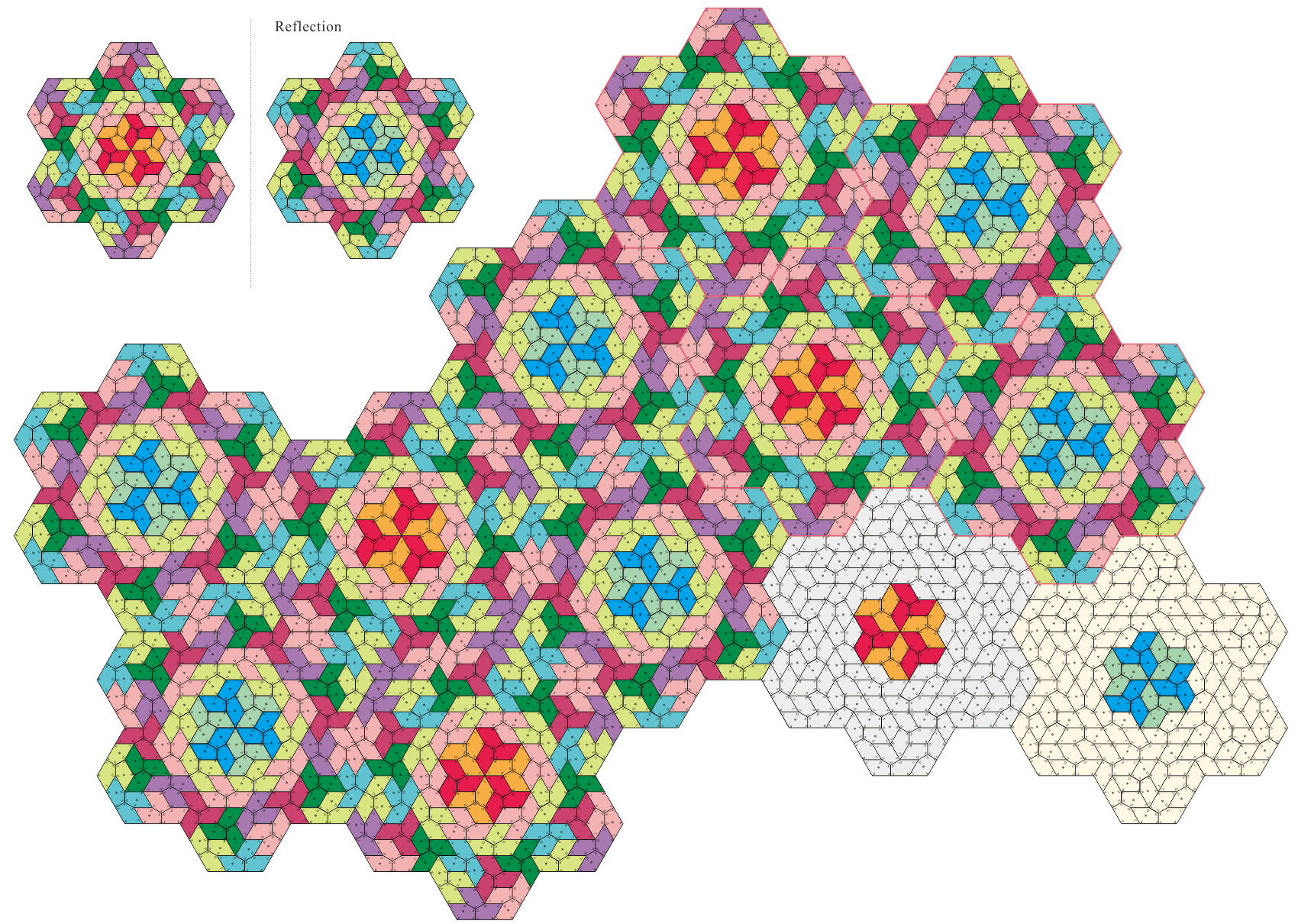} 
  \caption{{\small 
Example of tiling by HFL3-units.} 
\label{fig50}
}
\end{figure}

\begin{figure}[htbp]
 \centering\includegraphics[width=15cm,clip]{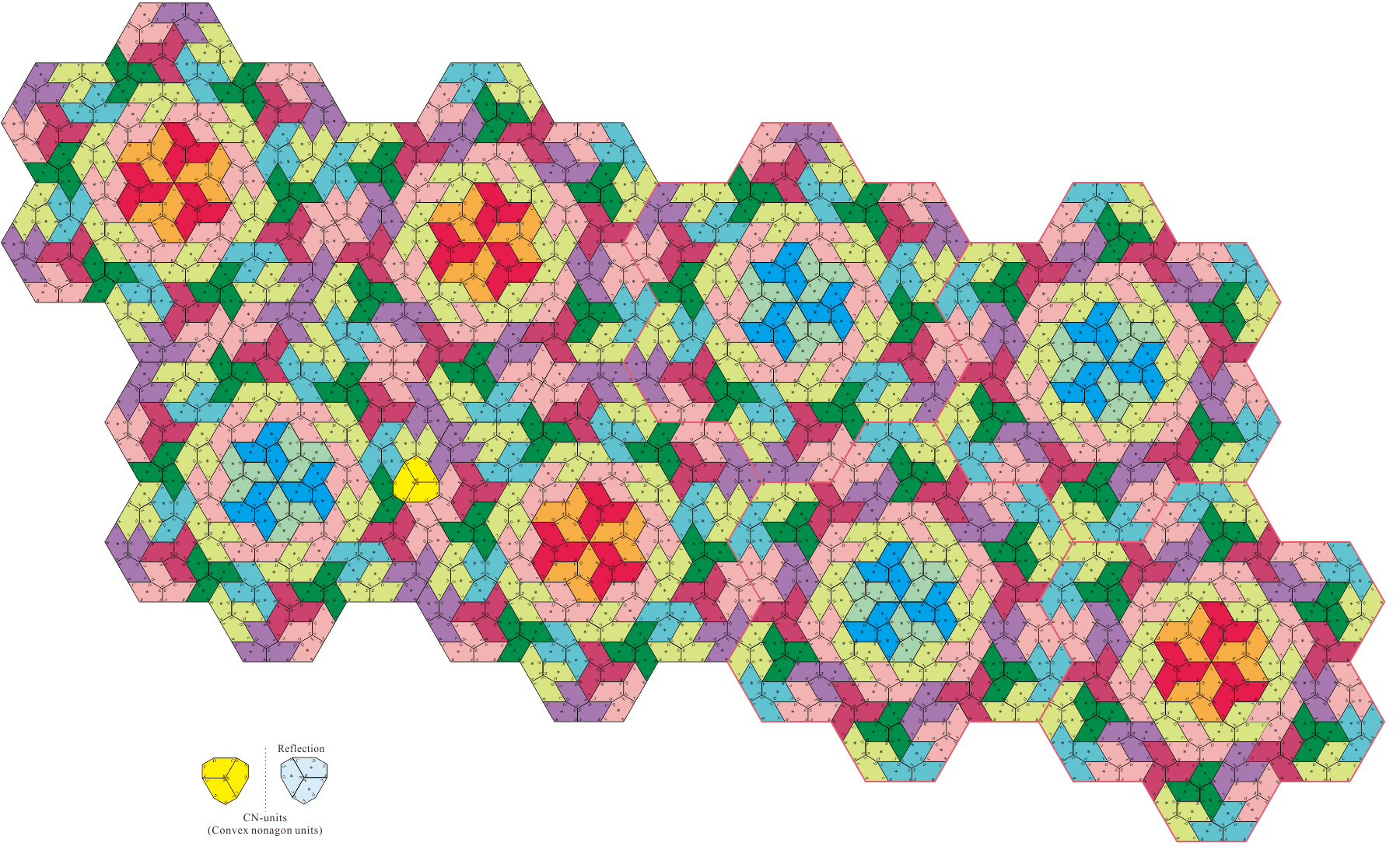} 
  \caption{{\small 
Example of tiling by HFL3-units.} 
\label{fig51}
}
\end{figure}

Like HFL1-units, the HFL3-units can generate tilings by combining 
HFL3-units. The contiguity methods are three patterns as shown in 
Figures~\ref{fig48}, \ref{fig49}, and \ref{fig50}. Note that the tilings in 
Figures~\ref{fig48}, \ref{fig49}, and \ref{fig50} are formed by an 
HFL3-unit of Figure~\ref{fig47}(a) and its reflected image. Then, the parts of 
CN-units in tilings can be reversed freely, and different patterns of 
HFL3-units which contain reflections can be freely connectable in the one 
tiling (see Figure~\ref{fig51}).

\renewcommand{\figurename}{{\small Figure.}}

\renewcommand{\figurename}{{\small Figure.}}

\renewcommand{\figurename}{{\small Figure.}}

\renewcommand{\figurename}{{\small Figure.}}

\renewcommand{\figurename}{{\small Figure.}}

\subsection{Skewed hexagonal flower unit and tilings}
\label{subsection5_5}

\begin{figure}[htbp]
 \centering\includegraphics[width=15cm,clip]{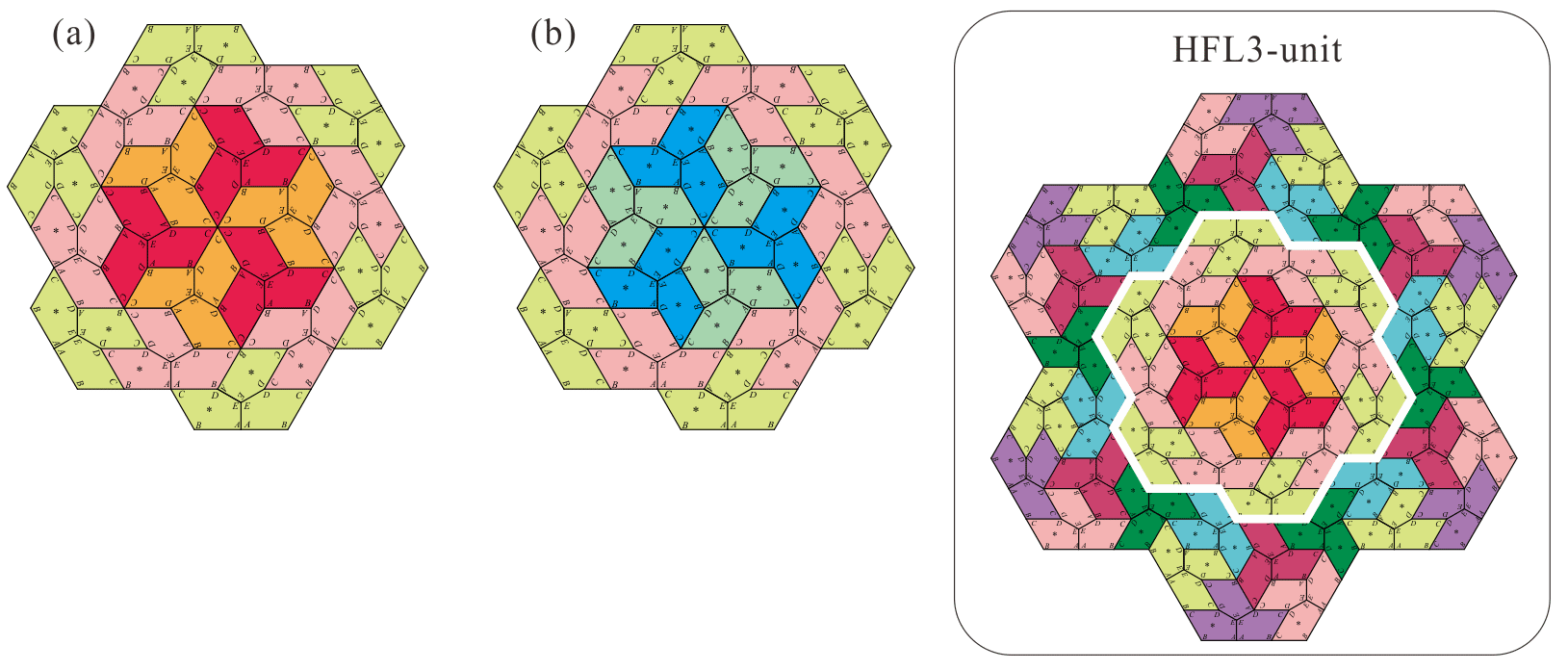} 
  \caption{{\small 
Skewed HF-units of two unique patterns.} 
\label{fig52}
}
\end{figure}

From the inside of the HFL3-unit, take out the 18 sided polygon as shown in 
Figure~\ref{fig52}. This 18-sided polygon is formed with side lengths of two and one. 
Therefore, this 18-sided polygon is called a skewed hexagonal flower unit 
(skewed HF-unit). The unique arrangement of the internal TH-pentagon (or 
windmill unit and ship unit) that forms a skewed HF-unit is of two types 
shown in Figure~\ref{fig52}. The skewed HF-unit that contains 54 TH-pentagons, is 
made of an HFL1-unit, six windmill units, and twelve ship units. 

\begin{figure}[htbp]
 \centering\includegraphics[width=15cm,clip]{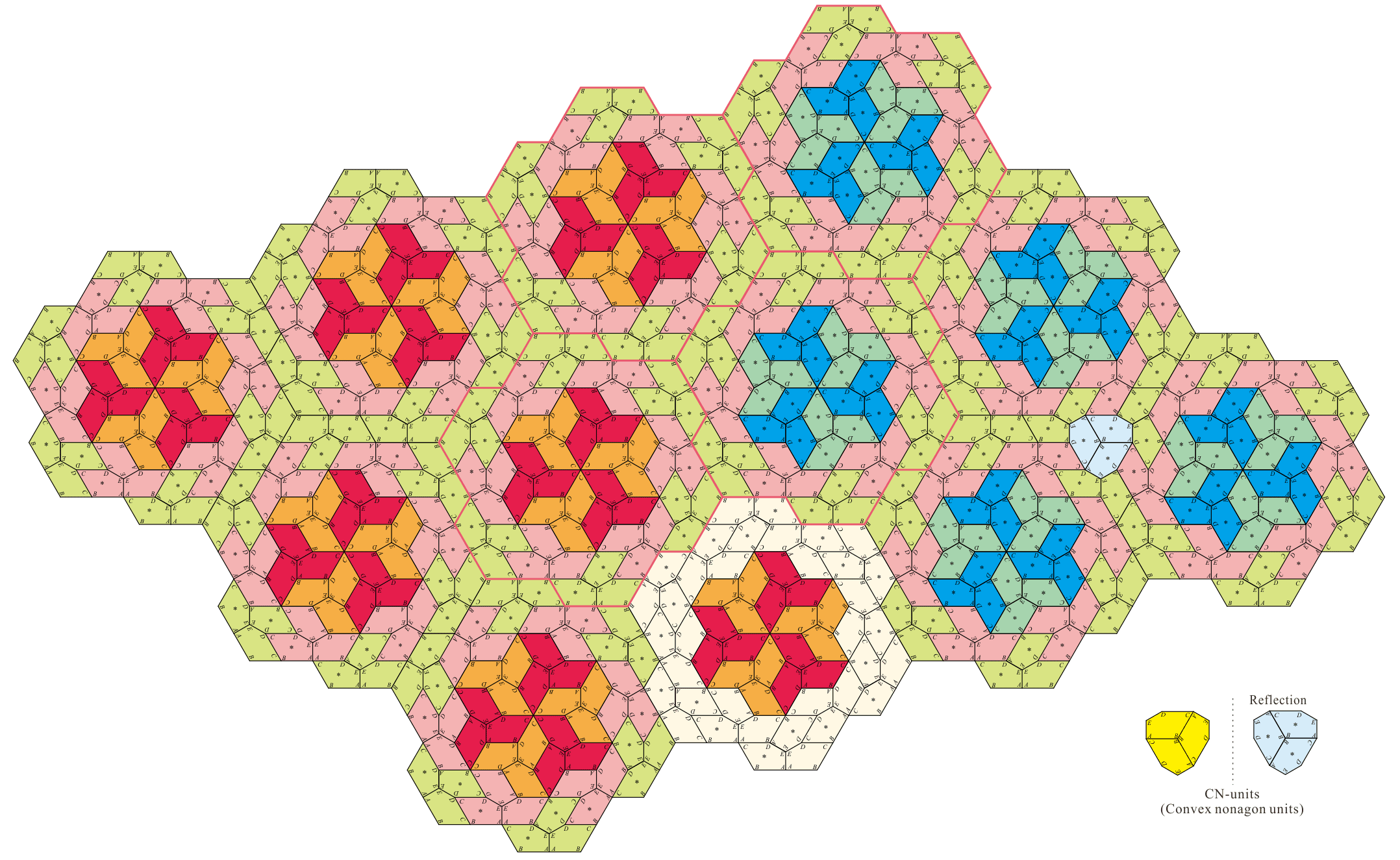} 
  \caption{{\small 
Example of tiling by skewed HF-units.} 
\label{fig53}
}
\end{figure}

The skewed HF-unit has only one contiguity method for a tiling. However, the 
parts of CN-units in tilings can be reversed freely, and different patterns 
of skewed HF-units can be freely connectable in the one tiling (see Figure~\ref{fig53}).

Using a program, it was confirmed that there is no TH-pentagon arrangement 
that fills a similar 18 sided polygon in Figure~\ref{fig52}.

\renewcommand{\figurename}{{\small Figure.}}

\renewcommand{\figurename}{{\small Figure.}}

\subsection{Tiling by the windmill units and the ship units based on the tiling 
with Classes S2 and S4}
\label{subsection5_6}

\begin{figure}[htbp]
 \centering\includegraphics[width=15cm,clip]{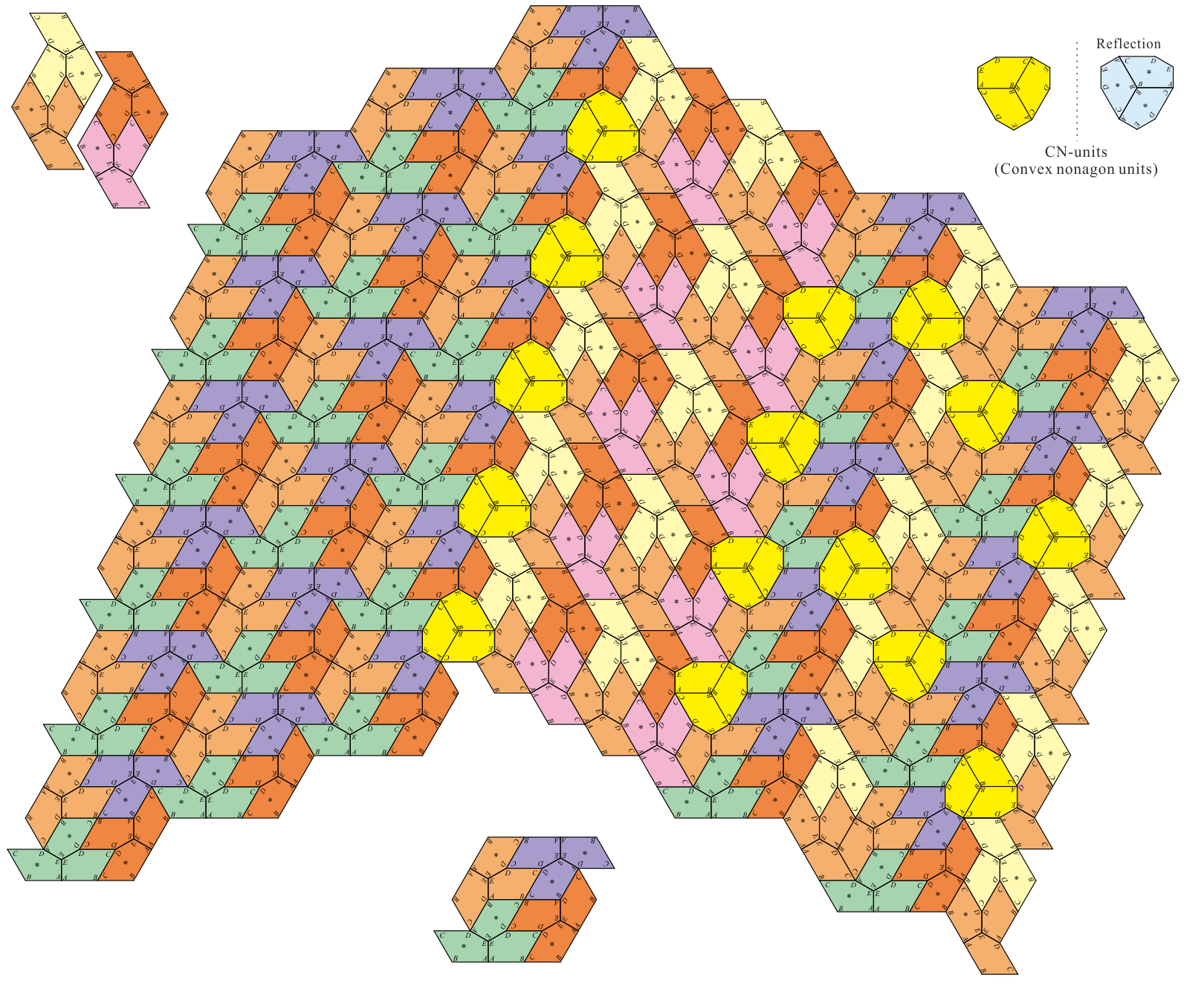} 
  \caption{{\small 
Tiling by the windmill units and the ship units based on the 
tiling with Classes S2 and S4.} 
\label{fig54}
}
\end{figure}

The tiling with Classes S2 and S4 in Figure~\ref{fig38} in Section~\ref{section4} contains 
CN-units. Therefore, as shown in Figure~\ref{fig54}, since the parts of CN-units in 
tilings can be reversed freely, the tiling of Figure~\ref{fig38} by only the ship 
units can be turned into a tiling with windmill units and ship units.

\renewcommand{\figurename}{{\small Figure.}}

\section{Conclusion }

The authors identified novel properties of the tilings of TH-pentagon. As a 
result, many new tilings were found (see Sections~\ref{section3}, ~\ref{section4}, 
and ~\ref{section5}). The TH-pentagon admits many periodic tilings and 
nonperiodic tilings. Thus, the TH-pentagon can create infinite tilings.

Since the search has not ended, new tilings will be found in the 
future\footnote{ After writing this manuscript, it was found that Johannes 
Hindriks presents tilings by heptiamonds at the site 
\url{http://www.jhhindriks.info/37/}. The authors confirm that there are tilings 
which they have not found yet. They are introduced in the following 
manuscript.}. If a polyhex ($n$-hexes) based on hexagons can be filled by a 
windmill unit and a ship unit, there may be a new convex pentagon 
tiling~\cite{G_and_S_1987, wiki-15, wiki-16}. The shape of a hexagonal 
flower unit corresponds to one of Heptahexa (7-hexes).

\bigskip

\bigskip
\noindent
\textbf{Acknowledgments.} 
The authors would like to thank Mr. Toshihiro Shirakawa for helpful comments.

\appendix
\def\thesection{Appendix }
\section{}

\renewcommand{\figurename}{{\small Figure.}}
\begin{figure}[htbp]
 \centering\includegraphics[width=15cm,clip]{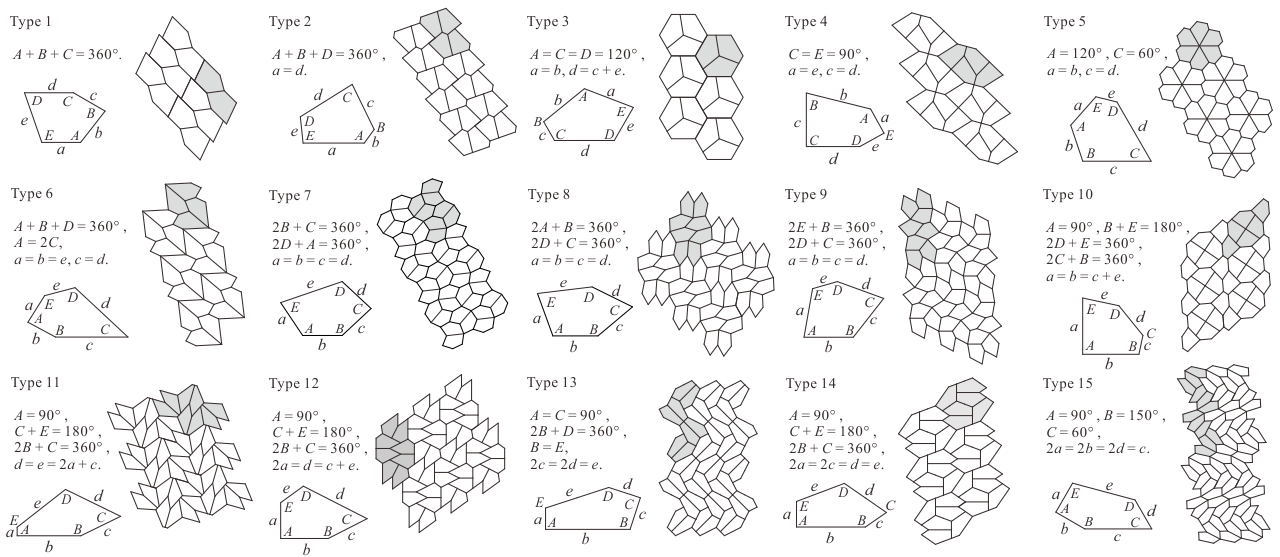} 
  \caption{{\small 
 Convex pentagonal tiles of 15 types. Each of the convex 
pentagonal tiles is defined by some conditions between the lengths of the 
edges and the magnitudes of the angles, but some degrees of freedom remain. 
For example, a convex pentagonal tile belonging to Type 1 satisfies that the 
sum of three consecutive angles is equal to $360^ \circ$. This condition for 
Type 1 is expressed as $A+B+C=360^ \circ$ in this figure. The pentagonal 
tiles of Types 14 and 15 have one degree of freedom, that of size. For 
example, the value of $C$ of the pentagonal tile of Type 14 is 
$\cos ^{ - 1}((3\sqrt {57} - 17) / 16) \approx 1.2099\;$rad $\approx 69.32^ \circ $.} 
\label{fig55}
}
\end{figure}

\renewcommand{\figurename}{{\small Figure.}}
\begin{figure}[htbp]
 \centering\includegraphics[width=15cm,clip]{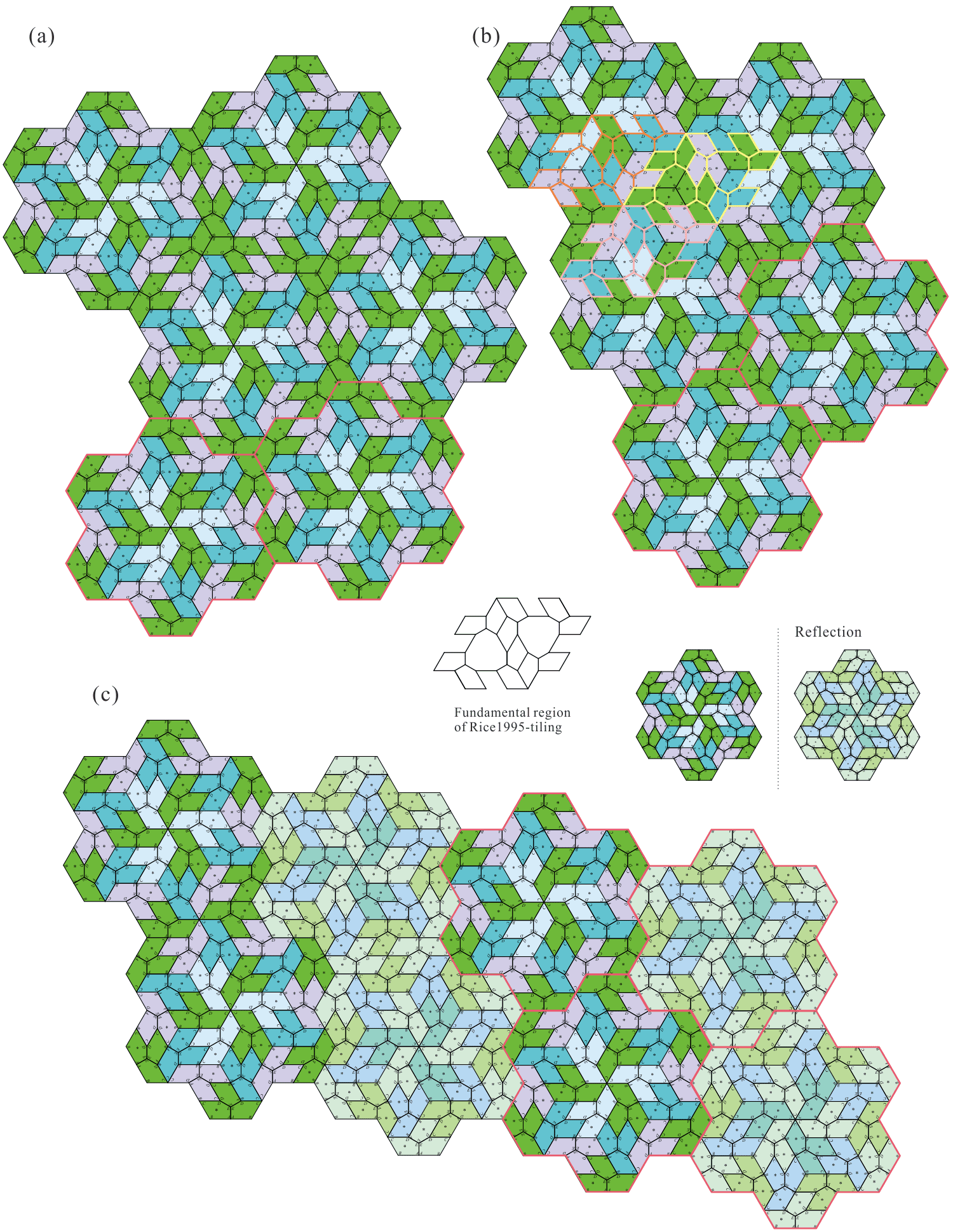} 
  \caption{{\small 
Tilings by HFL2-unit of Figure ~\ref{fig41}(n) and its reflection image.} 
\label{fig56}
}
\end{figure}

\renewcommand{\figurename}{{\small Figure.}}
\begin{figure}[htbp]
 \centering\includegraphics[width=15cm,clip]{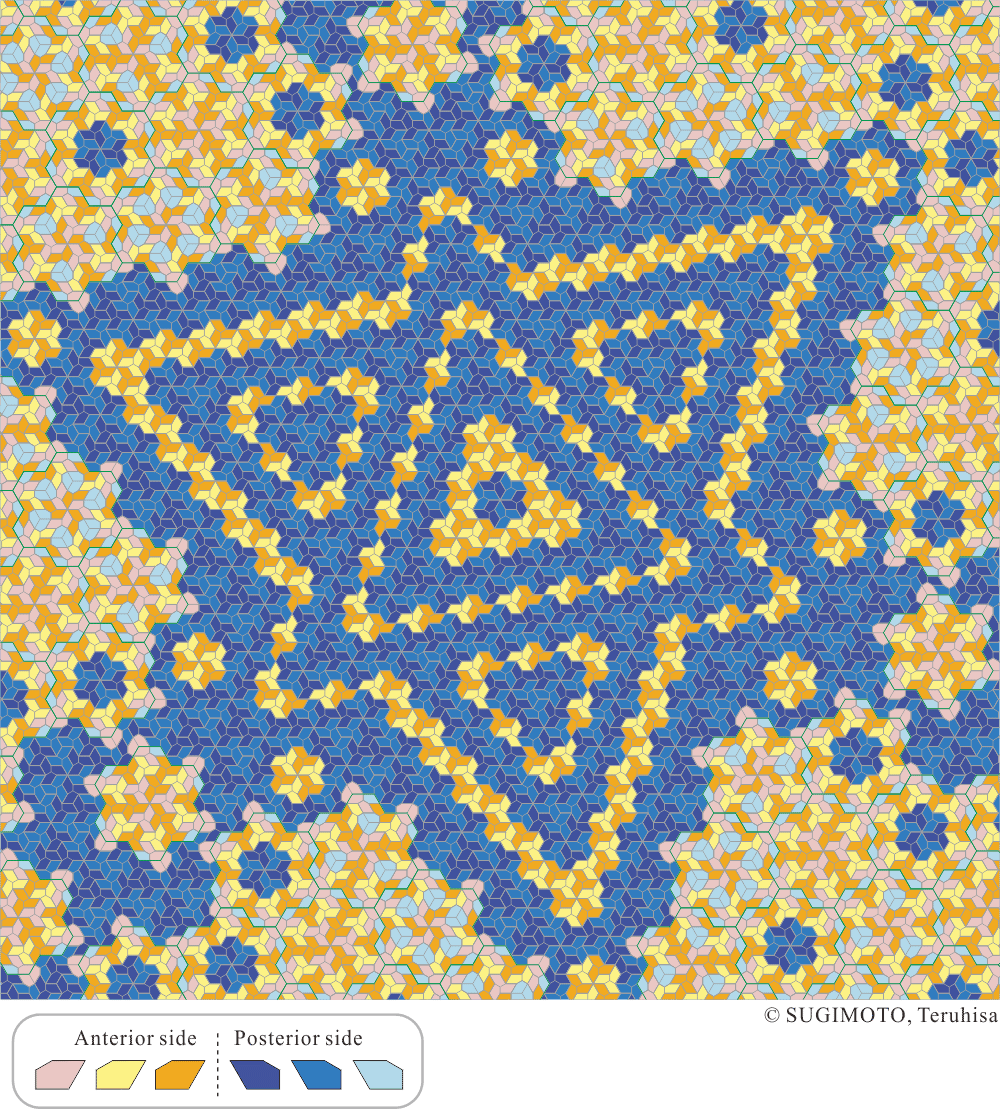} 
  \caption{{\small 
Tiling by windmill units, HFL1-units, and HFL2-units.} 
\label{fig57}
}
\end{figure}

\renewcommand{\figurename}{{\small Figure.}}
\begin{figure}[htbp]
 \centering\includegraphics[width=15cm,clip]{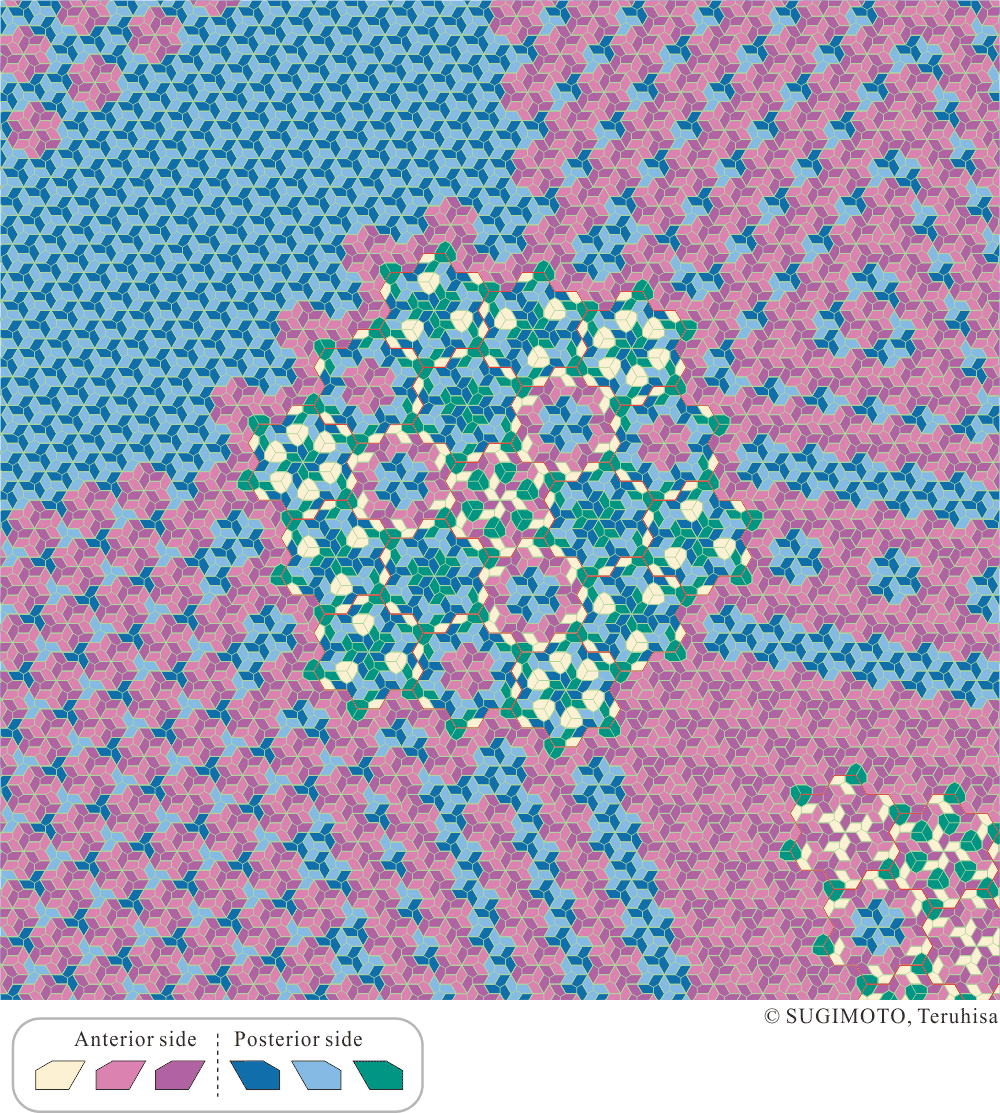} 
  \caption{{\small 
Tiling by windmill units, HFL1-units, and HFL2-units}
\label{fig58}
}
\end{figure}

\end{document}